\documentclass[11pt]{amsart}
\usepackage{graphicx}
\usepackage{cite,amssymb,amsmath}
\usepackage{amsthm}

\newtheorem{theorem}{Theorem}[section]
\newtheorem{corollary}[theorem]{Corollary}
\newtheorem{lemma}[theorem]{Lemma}
\newtheorem{proposition}[theorem]{Proposition}

\theoremstyle{definition}

\theoremstyle{remark}

\newcommand{\newoperator}[2]{\DeclareMathOperator{#1}{#2}}
\newcommand{\approxspace}[1]{{}}

\newcommand{\df}[1]{\emph{#1}}
\newoperator{\AH}{AH}
\newoperator{\AJ}{AJ}
\newoperator{\BS}{BS}               % branched surface
\newoperator{\GF}{GF}               % Hausdorff quotient of branched surface.
\newoperator{\GH}{GH}
\newoperator{\GJ}{GJ}
\newoperator{\Isom}{Isom}
\newoperator{\KS}{KS}               % Hausdorff quotient of branched surface.
\newoperator{\Nbhd}{Nbhd}
\newoperator{\PL}{PL}
\newoperator{\QF}{QF}
\newoperator{\QH}{QH}
\newoperator{\QJ}{QJ}
\newoperator{\QpH}{Q^\prime H}
\newoperator{\Rad}{Rad}
\newoperator{\Tr}{Tr}
\newoperator{\Vol}{Vol}
\newoperator{\area}{area}
\newoperator{\base}{base}
\newoperator{\codim}{codim}
\newoperator{\const}{const}
\newoperator{\degree}{degree}       % as in degree of a map
\newoperator{\glue}{glue}
\newoperator{\graft}{graft}
\newoperator{\inj}{inj}             % injectivity radius
\newoperator{\interior}{int}
\newoperator{\length}{length}
\newoperator{\llength}{\underline length}
\newoperator{\ls}{ls}
\newoperator{\qf}{qf}
\newoperator{\skin}{skin}
\newoperator{\teich}{\mathcal{T}}
\newoperator{\unglue}{unglue}
\newoperator{\wb}{wb}               % window base: surface is 2-fold covering.
\newoperator{\window}{window}
\newoperator{\ws}{ws}               % window surface

\newcommand{\hy}{\mathbb{H}}

\newcommand{\complexes}{\mathbb{C}}

\newcommand{\reals}{\mathbb{R}}
\newcommand{\euclidean}{\mathbb{E}}

\newcommand{\arrow}{\rightarrow}

\newcommand{\cross}{\times}

\newcommand{\homotopic}{\simeq}

\newcommand{\inverse}{^{-1}}
\newcommand{\product}{\prod}

\newcommand{\union}{\cup}

\DeclareMathOperator{\Stretch}{{ stretch}}
\DeclareMathOperator{\far}{{far}}
\DeclareMathOperator{\cut}{{cut}}
\DeclareMathOperator{\cat}{{cat}}
\DeclareMathOperator{\Cut}{{Cut}}
\DeclareMathOperator{\Near}{{Near}}
\DeclareMathOperator{\hol}{ {hol}}
\DeclareMathOperator{\CL}{\mathcal{CL}}
\DeclareMathOperator{\ML}{\mathcal{ML}}
\DeclareMathOperator{\MF}{\mathcal{MF}}
\DeclareMathOperator{\PF}{\mathcal{PF}}
\newcommand{\IP}{\par \noindent \  $\bullet$ \quad}
\newcommand{\element}{\in}
\begin{document}
\title{ Minimal stretch maps between hyperbolic surfaces }
\author {William P. Thurston}
\address{
Mathematics Department \\
University of California at Davis \\
Davis, CA 95616
}
\email{wpt@math.ucdavis.edu}
\begin{abstract}
Teichm\"uller maps between surfaces have had significant 
mathematical applications. This paper develops a geometric theory
that parallels the theory of Teichm\"uller maps. The central result is that for
any two hyperbolic structures on a surface of finite area, the
minimum Lipschitz constant for a map between them equals the supremum
(which with probability one is the maximum)
of the ratio of lengths of simple closed geodesics in the two
surfaces. Associated to this phenomenon is an asymmetric Finsler metric
on Teichm\"uller space. The extremal Lipschitz maps between surfaces are
closely related to geodesics in this metric; both geodesics and maps
are constructed.

In the process, some geometrically-defined \df{cataclysm} coordinate systems
are developed for Teichm\"uller space and measured lamination spaces,
along with properties and formulas for differentiation. 
These coordinates and computations have much wider applicability than
to the particular application of analyzing extremal Lipschitz maps.
Teichm\"uller space is described as the interior of a convex ball, in a way that
the boundary is its projective measured lamination compactification,
and so that the action of the modular group is differentiable
in a suitable sense.

This eprint is a 1986 preprint converted to \LaTeX with only
minor corrections. I currently think that a characterization
of minimal stretch maps should be possible in a considerably more general
context (in particular, to include some version for all Riemannian surfaces),
and it should be feasible with a simpler proof based more
on general principles---in particular, the max flow min cut principle,
convexity, and $L^0 \leftrightarrow L^\infty$ duality.
However, the explicit finite-dimensional nature of the present paper
seems interesting for its own sake, so I have refrained from working on
any significant revision.  I also suspect that this theory fits into a context
including $L^p$ comparisons, and that there should be complexified versions
connecting to the study of isometry classes of immersions of surfaces
in hyperbolic manifolds, and that it has connections to the quasi-isometric
study of leaves of foliations, pseudo-Anosov flows, \emph{etc.} (see
\cite{Thurston:circlesI} \dots .)

Some basic background is developed in \cite{Thurston:3dgt} and \cite{gt3m}.
\end{abstract}
\date{1986 preprint converted to 1998 eprint}
\maketitle
\sloppy
\tableofcontents
%\listoffigures
\section {Introduction} \label {Introduction}

A topological surface can have many different geometric shapes.
These variable shapes are recorded and measured by
the Teichm\"uller space of the surface,
which can be defined
either as the set of hyperbolic structures on a surface, or as
the set of conformal structures on a surface.

Teichm\"uller space is more than just a set.
It has a topology, which can
be defined in a variety of ways.
A neighborhood basis for an element $x$
of Teichm\"uller space might be
defined by formalizing any of the following descriptions:
\begin{itemize}
\item
conformal structures related to $x$
by a diffeomeorphism
isotopic to the identity with small conformal distortion (dilatation).
\item
conformal structures such that for a finite set of simple closed curves
the extremal lengths
of those curves are near the extremal lengths for $x$.
\item
hyperbolic structures for which the lengths of a finite set of
geodesics are nearly the same as in $x$.
\item
hyperbolic structures such that
no length of any geodesic has changed by more than a small ratio.
\item
hyperbolic structures such that
no length of any geodesic has increased by more than a small ratio.
\item
hyperbolic structures
for which there exists a diffeomorphism whose derivative
has $L^{p}$ norm near 1.
\end{itemize}

Despite the variety of possible reasonable definitions --- the list
goes on and on ---
the resulting topology is the same.
Topologically, Teichm\"uller space is the same as
Euclidean space of the appropriate dimension.

Teichm\"uller space is also more than just a topological space.
In order to truly understand the totality of possible geometric
structures on a surface,
we would like a geometric understanding of Teichm\"uller space itself.
When we begin to try to impose a geometric structure on Teichm\"uller space,
the different possible definitions
begin to diverge not just in their flavor, but in their results.
Teichm\"uller space comes equipped not with one geometry, but with
a number of distinct geometries. %(See REFS).
These different geometries seem to be associated with different kinds of
questions that arise concerning the geometry of surfaces.
Experience indicates that no one of these geometries is sufficient:
different geometries are helpful in different contexts.
Perhaps some day we will understand how all the different structures on
Teichm\"uller space fit together, but that day has not yet arrived.

Some of the geometries have been thoroughly studied, and we are
beginning to have a fairly clear picture of them.
For example, the Teichm\"uller metric arises from the question
``What is the least possible distortion of conformal structure for a
quasiconformal map between two surfaces with a fixed homotopy class?''
Teichm\"uller proved that there is a unique quasiconformal map with
minimum possible quasiconformal distortion.
This map is known as the Teichm\"uller map.
The amount of distortion gives rise to a metric, the Teichm\"uller metric,
which has geodesics that are described by one-parameter families of
Teichm\"uller maps.
The Teichm\"uller metric turns out to be the same as the Kobayashi
metric, which is defined for a large class of complex manifolds
\cite {Royden}.

Another classical metric is the Weil-Peterson metric, a Riemannian metric
on Teichm\"uller space
which is related to harmonic maps between surfaces (see for instance
\cite {Wolf:harmonic}, \cite {Wolpert2}.)
Through the efforts of Scott Wolpert and others, we now have a fairly clear
picture of this metric.
It is a metric of strictly negative sectional curvature.
It is not a complete metric: a family of hyperbolic surfaces
in which a simple closed geodesic shrinks to nothing (in a direct way)
represents a path in Teichm\"uller space
which has finite length as measured by the Weil-Peterson metric.
All is not lost, however, because it turns out that
Teichm\"uller space is geodesically convex
with respect to the Weil-Peterson metric.
This means that any two points are joined by a unique geodesic, and
many global geometric constructions are still possible \cite{Wolpert6}. 

In this paper, we will investigate a different kind of
geometry for Teichm\"uller space
which involves the lengths of geodesic laminations on a surface and the
Lipschitz norm (or 
$L^{\infty}$
norm on the derivatives) for
maps between surfaces.
In particular, we will study a
question reminiscent of the question which was
answered by Teichm\"uller:
\IP
Given any two hyperbolic surfaces 
$S$ and $T$,
what is the least possible value of the global Lipschitz constant
$$
L(\phi ) = \sup_{x \ne y} \left ( \frac{ d( \phi (x), \phi (y)) }{ d(x, y) }
\right ) $$
for a homeomorphism 
$\phi :  S \arrow T$
in a given homotopy class?

This is closely related to the canonical
problem that arises when a person on the standard American diet
digs into his or her wardrobe of a few years earlier.
The difference is that in the wardrobe problem,
one does not really care to know
the value of the best Lipschitz constant --- one is mainly concerned that the
Lipschitz constant not be significantly greater than 1.
We shall see that,
just as cloth which is stretched tight develops stress wrinkles,
the least Lipschitz constant for a homeomorphism
between two surfaces is dictated by a certain geodesic lamination
which is maximally stretched.

For any $K$
the set of $K$-Lipschitz maps between two compact metric spaces
is equicontinuous, and hence compact.  It follows
that there is a continuous map from $S$ to $T$
which realizes the best constant, but it is
not so easy to see what geometric properties it might have.
In particular, it is not
obvious that it can be taken to be a homeomorphism.

We will give a geometric construction for a homeomorphism
between the two hyperbolic surfaces which does have the best
Lipschitz constant.  The Lipschitz-extremal
homeomorphism is not as nice as the Teichm\"uller map --- in particular,
it is not unique --- but still it has a geometric form which is
strongly reminiscent of the Teichm\"uller map.

If $M^{2}$ is a surface and if $g$ and $h$
are complete hyperbolic structures on $M^{2}$
representing two elements of Teichm\"uller space $\teich(M^2)$, then
the Lipschitz comparison between the two metrics is conveniently measured
by the quantity
$$
L ( g , h ) = \inf_{\phi \homotopic \text {id} } \left ( \log ( L ( \phi ) ) \right )
$$
where it is understood that in computing the Lipschitz constant 
$L ( \phi )$,
the metric $g$ is used before $\phi$ and the metric 
$h$ is used after $\phi$.

The construction for a Lipschitz-extremal homeomorphism involves the analysis of
another measure of the difference between the geometry of
two hyperbolic structures.  If $\alpha$ is any homotopically non-trivial
loop on $M^{2}$, it is represented by a
geodesic in each of the two metrics.  Define a quantity
$$ K(g,h) = \sup_\alpha \left ( \log\left (
\frac{ \length ( \alpha_g ) }{\length ( \alpha_h ) }\right ) \right )
$$
Obviously, we have
$$ K(g,h) \le L(g,h) 
$$
since a Lipschitz map with Lipschitz constant 
$A$ sends a geodesic
to a curve whose length is multiplied by not more than $A$.
We shall prove that the two quantities are in fact equal.
This is a nice situation.
Even without knowing that the two quantities are equal,
the inequality gives us a concrete
way to bound either one: from above by producing a
map with a certain Lipschitz constant,
or from below by producing a geodesic whose ratio of lengths
is a certain amount.

There is a straightforward extension of the length function on
closed geodesics to a length function for geodesic laminations.
The length of a geodesic lamination is the total mass of the
product of its transverse measure with 1-dimensional Lebesgue measure
along the leaves.
Length is a continuous function on the space of measured geodesic
laminations.  Since the set of measured geodesic laminations, up to scalar
multiples, is compact, $K(g,h)$
is more concretely realized as
the log of the ratio of lengths of some measured geodesic
lamination in the two metrics.

One complication in the theory stems from the fact that the measured
geodesic lamination which attains the maximal ratio is not unique in
general, even up to scalar multiples.
We will prove a weaker uniqueness property, however:
there is a topological lamination $\mu (g, h)$
so that the set of measured laminations attaining $K(g,h)$
consists of exactly those measured laminations defined by transverse
invariant measures (not necessarily of full support) on $\mu ( g, h)$.

In the course of this paper, we will develop some other ideas which are
of interest in their own right.  The intent is not to give the slickest
proof of the main theorem, but to develop a good picture.
We will describe a model for measured lamination space as the boundary
of a certain convex cone.
We will define and analyze the tangent space of measured lamination space.
We will give a geometric calculation for the derivative of the length
function on the space of measured laminations, insofar as there is
a derivative.
We will also introduce a generalization of the notion of an earthquake
on a surface, namely the cataclysm.
Among other things, cataclysms give rise to global coordinate systems
for Teichm\"uller space which have nice geometric properties.

\section {Elementary properties of the Lipschitz constant} \label {Elementary properties of the Lipschitz constant}

Let $S$ be a surface, and $\teich(S)$ its Teichm\"uller space.
The function 
$L(g,h) $ on $\teich(S) \cross \teich(S)$
can be thought of as one way of
measuring the distance between two metrics.

\begin{figure}[htbp] \centering \approxspace{ 4in} 
\includegraphics{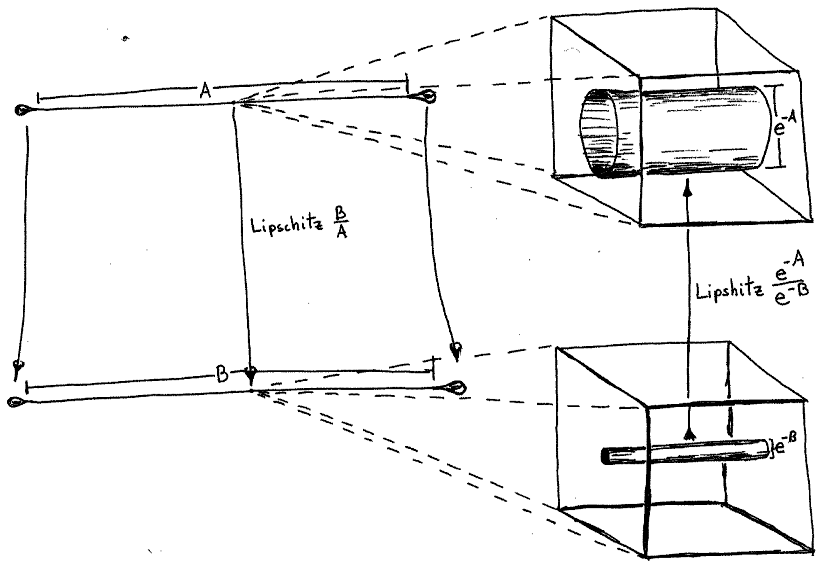}
\caption{
\label{asymmetry of L} 
The best Lipschitz constant for a map between two surfaces is
not a symmetric invariant of the pair of surfaces.
Illustrated on the left are two surfaces of genus 2; seen from
a distance, their handles appear
quite small, like eyes of needles. The ratio of lengths of the
two-eyed needles is approximately 1, but the ratio of their diameters
is large. These two ratios are the optimal Lipschitz constants for maps in
the two directions.}
\end{figure}

This distance function is not actually a metric, since
it is not symmetric: it rarely happens that 
$L(g,h) = L(h,g)$.
The asymmetry is vividly illustrated by the
example (figure \ref{asymmetry of L})
of a hyperbolic surface with a very short geodesic.
Compare it to another surface which is very nearly the same in the
thick portion, but where the length of the tube about the short geodesic
has been increased.  If there is a tube of length $A$
about the short geodesic in the first metric and length 
$B$ in the second metric, there is a map from the
first to the second with Lipschitz constant approximately
$$ \frac BA $$
However, the lengths of the short geodesics have ratio approximately
$$ \frac {\exp ( - B ) }{ \exp ( - A ) } $$
so that any map in the reverse direction must have a huge
Lipschitz constant.

The other properties of a metric are satisfied by $L$.
It is immediate from the definition
that it satisfies the triangle inequality
$$
L(f,g) + L(g,h) \ge L(f,h).
$$
The order of the arguments matters, because
of the asymmetry.
This inequality is the reason for using
logs in the definition.  Finally,
\begin{proposition}[L distinguishes points]\label{L distinguishes points}
For any two hyperbolic structures,
$$ L(g,h) \ge 0 , $$
and
$$ L(g, h) = 0 \quad \text{only when} \quad g = h .$$
\end{proposition}
\proof %{Proof of \shortlabel L distinguishes points:2.2:}
For any $g,h$ such that $L=L(g,h) \le 0$
we can pick a continuous map with minimal global Lipschitz
constant $L$.
Every disk of small diameter in the domain metric is mapped to a subset of a
disk of the same diameter in the range surface.
But the area of the two surfaces is the same.
Covering the domain by a disjoint union of disks of full measure,
one sees that the only way for this to happen is that each disk
maps surjectively to a disk of the same size, so the map is an isometry.
\endproof %{\ref {L distinguishes points}}

It would be easy to replace $L$ by its symmetrization
$ 1/2 ( L(g,h) + L(h,g))$, but it seems that, because of its direct
geometric interpretations, $L$ is more useful just as it is.

One consequence of our construction will be that there is a non-symmetric
Finsler metric, that is, a non-symmetric norm on the tangent bundle of
$\teich(S)$, such that $K$ is minimal path length for this norm.

There are other quantities one might try to minimize,
to get other kinds of extremal maps.
Perhaps the most natural question in this vein is whether there is
a good description for a best bi-Lipschitz homeomorphism, that is,
a homeomorphism which minimizes the
larger of its Lipschitz constant and the Lipschitz constant of its inverse.
The existence of such a homeomorphism is immediate.

It is worth mentioning here that the construction we will make for a
Lipschitz-extremal homeomorphism between two surfaces has an inverse
which is not Lipschitz.  As the discussion proceeds, it will become
evident that for most pairs of surfaces it is not possible for
a Lipschitz-extremal map to have a Lipschitz inverse.
This is already visible in the basic ingredient we shall use
for constructing Lipschitz maps:
\begin{proposition}[Triangle expands]\label{Triangle expands}
For any constant $K > 1$ there is a $K$-Lipschitz
homeomorphism of a (filled) ideal hyperbolic triangle to itself
which maps each side to itself, and multiplies arc length on the sides
by the constant $K$.
\end{proposition}

\begin{figure}[htbp] \centering \approxspace{ 4in} 
\includegraphics{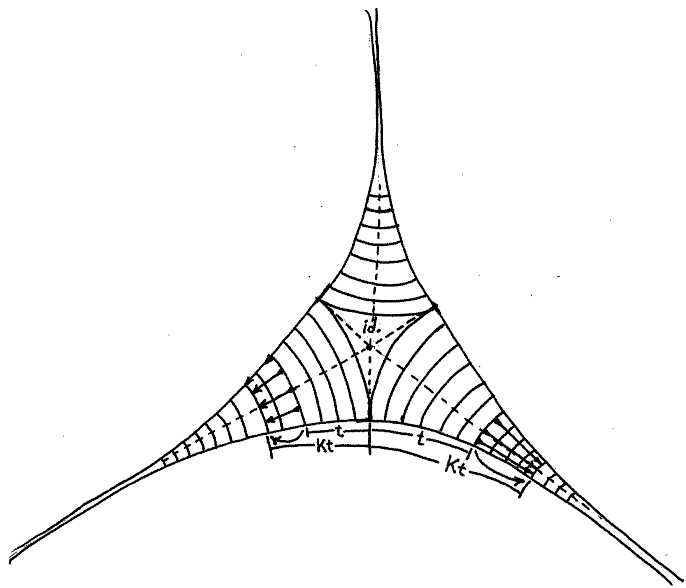}
\label{expansion of  triangle} 
\caption{ 
For each $K > 1$ there is a $K$-Lipschitz map of an ideal triangle
to itself which exactly expands each side by 
$K$.  }
\end{figure}

\proof %{Proof of \shortlabel Triangle expands:2.3:}
Each corner of an ideal triangle can be foliated by horocycles.
Extend these foliations until they fill all but
region in the center bounded by three horocycles.

A $K$-Lipschitz homeomorphism can be constructed by
fixing the central region, and mapping
a horocycle which has distance $t$ from the central region to
the horocycle with distance $K t$.
The map can be made linear with respect to arc length on the horocycles.
\endproof %{\ref {Triangle expands}}

Note that the arc length of a horocycle is contracted by the
factor $\exp ( t ( K - 1 ))$, which becomes extremely
severe near the corners of the triangles.
If desired, the homeomorphism above could be relaxed in the interior of
the triangle so that it is differentiable, and the derivative has norm
strictly less than $K$ in the interior of the triangle.

\section {Properties of ratios of lengths} \label {Properties of ratios of lengths}

Given a pair of hyperbolic metrics $g$ and $h$ on a surface $S$,
there is associated with any element $\gamma  \element  \pi_{1} ( S ) $
a number
$$ r_{gh} ( \gamma ) = \frac { \length_h ( \gamma ) }{ \length_g ( \gamma ). }
$$
The quantity $K(g,h)$ which we defined earlier is the supremum of the
log of $r_{gh} ( \gamma )$, where 
$\gamma$ ranges over $\pi_{1} ( S )$.

We need to prove that $K ( g , h ) > 0$ for $g \ne h$.
One way to prove this would be to use the earthquake theorem, which asserts that
any two metrics $g$ and $h$ are related by a left earthquake.
Given such an earthquake, it is fairly easy to construct a geodesic which
whose length is increased going from $g$ to $h$, by following along nearly
tangent to the fault zone, but slowly crossing from left to right.

We will give a different proof, presenting a reason which
is more general and based more on first principles.

\begin{theorem}[K(g,h) positive]\label{K(g,h) positive}
For any two distinct hyperbolic structures $g$ and $h$,
$$ K ( g , h) > 0
$$
\end{theorem}
\proof %{Proof of \shortlabel K(g,h) positive:3.1:}

Let $g$ and $h$ be any two hyperbolic structures on $S$.
The two metrics have the same area, so there is an area-preserving
diffeomorphism $\phi : g \arrow h$ between them.
If there are any cusps, we can assume that $\phi$ is an isometry
outside some compact set.

Let $\tilde\phi : \hy^{2}\arrow \hy^{2}$ be a lift of $\phi$ to the universal
cover.  Choose a basepoint $p$ in the domain, and let $D$ be the
function on the hyperbolic plane defined as
$$ D(q) = d(\tilde \phi(q), \tilde \phi(p)) - d(q, p) .
$$

\begin{lemma}[D unbounded]\label{D unbounded}
If the function $D$ has an upper bound, then $g = h$.
\end{lemma}

\proof [Proof of Lemma \ref{D unbounded}]
It is well known that the map $\tilde\phi$ must extend to a homeomorphism
$\bar\phi$ defined on the compactified hyperbolic plane
$D^{2} = \hy^{2} \union S_{1}^{\infty}$.
A wedge between two rays through $p$ is sent by $\bar \phi$ to
an approximate wedge, bounded by two quasigeodesic rays which are
a bounded distance from two geodesic rays.  If $D$
has an upper bound, then the length of the arc
where the image of the approximate wedge hits the circle at infinity
cannot be less than
a multiple bounded above zero of the length of the arc of
the circle at infinity for the domain wedge:
the area wouldn't fit otherwise.

This says that the inverse of $\bar\phi$ is Lipschitz when restricted
to the circle at infinity.

We finish by an argument which was used by Mostow for the conclusion
of his rigidity theorems.  Since $\bar\phi \inverse$ is Lipschitz,
it has a derivative almost everywhere.
For almost any point on the circle at infinity, there is a sequence
of elements of the fundamental group of the circle whose derivatives
form an unbounded sequence of expansions.
The map $\bar\phi \inverse$ is invariant by conjugation in the
domain and range by the actions of the fundamental group of the surface.
If we carry out conjugation by the sequence, we see that 
$\bar\phi$ is arbitrarily close to a Moebius transformation,
therefore it is a Moebius transformation giving an isometry between $g$ and $h$.
\endproof %{\ref {D unbounded}}

\begin{lemma}[approximate orbits transitive]
\label{approximate orbits transitive}
Let $X$ be a volume-preserving vector field on a
manifold $M$ of finite volume which generates a flow for all time.
Then for any compact set $A \subset M$
and for any $\epsilon$ there is a constant $C$ such that any two points in $A$
are connected by an $\epsilon$-approximate flow line in time $< C$, where
an $\epsilon$-approximate flow line is a curve $\gamma(t)\  [0 < t]$
such that $ \| T\gamma(t) - X(\gamma(t))\| < \epsilon$.
\end{lemma}

\proof [Proof of Lemma \ref {approximate orbits transitive}]
For any $x \in A$, the set of points which can be reached from $x$ by an
$\epsilon$-approximate flow line of time $< T$ is an open set.  The Poincar\'e
recurrence lemma implies that if $T$ is large enough, this set includes $x$.

In fact, for each point $x$, there is a neighborhood of $x$ and a time $T$
such that any point in the neighborhood can be connected to any other
in time $< T$.  From compactness, it follows that there is some finite
open cover of $A$
by neighborhoods having this property, for some time $U$.  If there
are $N$ neighborhoods in the cover, then at least in time $NU$ there are
approximate flow lines from any point in $A$ to any other.
\endproof %{\ref {approximate orbits transitive}}

\proof [Continuation, proof of \ref{K(g,h) positive}]
To complete the argument, take any point $q$ in $\hy^{2}$
such that $D ( q )$ is quite large.  We may assume that $q$
maps to the thick part of $S$: since $\phi$
was chosen to be an isometry in a neighborhood of each cusp, any upper bound
in the unit tangent bundle over
a compact part of $S$ would carry over to an upper bound over all
of $S$.

Consider all extensions of the geodesic from $\tilde\phi ( p )$ to 
$\tilde\phi ( q) $ which have curvature less than some small constant 
$0 < \epsilon $.  These extensions come from approximate flow lines
of the geodesic flow, so by \ref{approximate orbits transitive}
there is an extension of bounded length which ends at $\tilde \phi(p)$.

We have constructed a closed curve for the hyperbolic structure $S$
with curvature less than $\epsilon$.  It is homotopic to a geodesic for $h$
which has about the same length.  Its preimage under $\phi$
is homotopic to a geodesic for $g$ which is shorter.
\endproof %{\ref {K(g,h) positive}}

We will prove that the supremum of $r_{gh}$ over simple closed curves is the
same as the supremum over the entire fundamental group.

First we study hyperbolic pairs of pants,
{\it i.e.}, hyperbolic structures with geodesic boundary on
the sphere minus three disks.

\begin{lemma}[Shrinking at the waist]\label{Shrinking at the waist}
If $g$ and $h$ are hyperbolic structures on a pair of pants, then
for any element $\alpha$ of its fundamental group which
is not freely homotopic to the boundary, if $r_{gh} ( \alpha ) > 1 $, it is
less than the maximum of $r_{gh}$ for boundary components.
\end{lemma}

\proof %{Proof of \shortlabel Shrinking at the waist:3.4:}
First, consider the case that $r_{gh} $ takes the same value $r$
for the three boundary components.
Decompose the domain into two hyperbolic triangles, with the spikes
spiralling around the three boundary components.

If the metric in each triangle is replaced by a metric induced by the map
of proposition \ref{Triangle expands}
so that each side is expanded by a factor $r$,
the resulting metrics piece together to form a new pair of
pants whose boundary components have been expanded by a factor
of $r$.  Since hyperbolic structures on a pair
of pants are classified by the lengths of the boundary
components, this determines a Lipschitz map to the target pair of pants.

Any interior geodesic is mapped to a curve of length less than $r$
times its original length, so this case is finished.

\medskip
For the general case, let $r$ be the maximum of $r_{gh}$
applied to the boundary components.  First construct an $r$-Lipschitz map to
the pair of pants whose boundary components are each
$r$ times as long as the original.
We will show that when one or more of its boundary components is
shrunk, interior lengths shrink as well.

To see this, extend the pair of pants to become a complete hyperbolic
surface (without boundary).  Let $\alpha$ be any infinite embedded
geodesic which intersects one of the boundary components
twice, so that both ends of $\alpha$ go out to infinity in the same end of the
surface.  If $\alpha'$ is a nearby geodesic, you can
cut out the strip between them and join
the resulting boundary components in such a way that endpoints of
their common perpendicular are identified.
The result is a surface whose convex core is a pair of pants
with two boundary components the same, and the third shorter by a
controllable amount.
There is a 1-Lipschitz map between the two complete hyperbolic surfaces.
(The infinite hyperbolic surface is necessary, otherwise the fact that
the two pairs of pants have the same area would block such a map).

A composition of these constructions demonstrates the proposition.

\endproof %{\ref {Shrinking at the waist}}
\begin{proposition}[supremum simple]\label{supremum simple}
The same supremum $K(g,h)$ for the log of $r_{gh} ( \alpha )$
is obtained if $\alpha$ only ranges over
homotopy classes represented by simple closed curves.
\end{proposition}

\proof %{Proof of \shortlabel supremum simple:3.5:}
Given any geodesic $\alpha$, if $\alpha$ is not simple we will
construct a shorter geodesic with a greater value of 
$r_{gh}$.  Clearly, that is sufficient to prove the theorem.

In the first place, we may assume that $\alpha$ is not a curve
which wraps more than once around its image.
Thus, all self-intersections are transverse.

Choose any self-intersection point $p$ of $\alpha$.
The geodesic describes an immersed figure 8, intersecting itself at $p$.
Form the covering space of the surface corresponding to the image of the
fundamental group of the 8.  The surface is a pair of pants, extended to
be a complete hyperbolic manifold.  Apply the preceding proposition to
complete the proof.
\endproof %{\ref {supremum simple}}

It is not hard to prove that if two simple closed curves $\alpha$ and 
$\beta$ both attain the supremum for $r_{gh}$,
then they cannot intersect, by a simple cut and paste argument.
However, it is only a special case when the supremum is
attained by a simple closed curve: we really need to prove the property for
laminations.  We shall do this in 
\S\ref {Maximal ratio laminations}, after we have developed the
theory of the tangent space to measured lamination space and derivatives
of lengths of laminations.  By then, it will be easy to understand how
well a lamination is approximated by simple closed curves.

Even though it happens only in a special case that the supremum of $r_{gh}$
is attained by a simple closed curve, we shall see in section
\ref{The maximum stretch lamination is almost always a curve}
that this special case is statistically the only significant one.

\section {Stretch maps between surfaces} \label {Stretch maps between surfaces}

In this section, we will describe the central construction, which
produces new hyperbolic structures
$$ \Stretch( g , \lambda , t )
$$
from a hyperbolic structure $g$ on a surface $S$
such that the Lipschitz constant for the identity map
and the maximum ratio of lengths of geodesics are in harmony.
The construction depends on

\IP a starting metric $g$.

\IP a geodesic lamination $\lambda$ whose complementary
regions are triangles, but does not necessarily carry a transverse invariant
measure of full support, and

\IP a real parameter $t$ which will satisfy
$$ 0 < t = K ( g , \Stretch ( g , \lambda , t )) =
L ( g , \Stretch ( g , \lambda , t )) .
$$

We assume that $g$ is a complete hyperbolic structure on $S$ of finite
area.  The lamination $\lambda$ necessarily has leaves going out to cusps
in the case that $S$ is not compact.

The construction amounts to the assembly of copies of the basic unit of
proposition \ref{Triangle expands}, one copy for each component of the
complement of $\lambda$.
This defines a change in the hyperbolic structure of each triangle
such that the sides of each triangle are expanded by a factor of 
$\exp ( t )$.  These changes match up along $\lambda$
so that they change the arc length of the leaves of $\lambda$
by a factor of $\exp ( t )$ as well.

\begin{figure}[htbp] \centering \approxspace{ 4in} 
\includegraphics{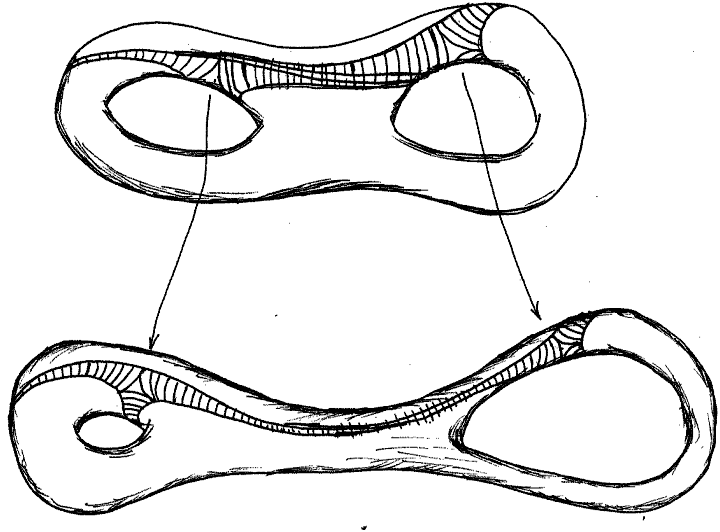}
\caption{  
Maps between certain pairs of hyperbolic surfaces can be constructed by
assembling the units of figure
in the complementary regions of a geodesic lamination.
}
\label{stretch map between surfaces} 
\end{figure}

The net result will be a new hyperbolic structure 
$h$, with a built in map from the original structure $g$ having Lipschitz
constant $\exp ( t )$.  If it is not the case that $\lambda$
consists entirely of leaves going out at both ends to cusps, then $\lambda$
admits a transverse invariant measure with compact support.
The corresponding measured lamination has its length multiplied by the
Lipschitz constant $\exp ( t )$, so no better Lipschitz constant is
possible.  Furthermore, we can approximate $\lambda$ in measured lamination
space by
a sequence of real multiples of simple closed curves; since length is
continuous in measured lamination space, it follows that
$K ( g , h ) = L ( g , h) = t$.

The remaining case, when all the leaves of $\lambda$ go out at both ends to
cusps, is the one case where the new metric $h$ is trivial to construct.
Unfortunately, this is also the one case where the built-in map between the two
metrics is not Lipschitz extremal: the map has a lot of unnecessary pull
toward the cusps.  This is clearly illustrated in the example of the
three-punctured sphere, where stretch maps do not change the underlying
metric at all since there is only one metric up to isometry.
More generally, whenever there are no leaves with compact closure,
it is possible to make an isotopy of the map which pushes horocycles
concentric to the cusps away from the cusps by a distance increasing linearly
toward the cusps.  This allows the map to be relaxed somewhat
everywhere, yielding a new map with lower Lipschitz constant.

In the general case, the intuitive explanation above for the construction of
$\Stretch $ is inadequate technically, because it does
not describe a hyperbolic structure in a neighborhood of $\lambda$.
We will next analyze the general form of a small hyperbolic neighborhood of a
geodesic lamination.  This analysis will tie in later with a construction for
more general transformations of hyperbolic structures, namely cataclysms,
which include both stretch maps and earthquake maps as special cases.

Consider, then, a geodesic lamination $\mu$ on a complete hyperbolic surface
of finite area.  We make assumptions neither about the existence of a transverse
invariant measure for $\mu$, nor about the topological or combinatorial
type of the complement of $\mu$.

Choose an $\epsilon$ small enough that the topological type of the 
$\epsilon$-neighborhood $N_{\epsilon} ( \mu )$ has stabilized.

A measured foliation $F ( \mu )$ can be constructed in $N_{\epsilon}$
so that the leaves of $F ( \mu )$ meet the leaves of $\mu$ orthogonally,
and so that the transverse measure for $F ( \mu )$ agrees with arc length
for $\mu$.  First construct the leaves of $F ( \mu )$ to be
horocycles in a neighborhood of each ideal vertex of $S$ cut by $\mu$.

\begin{figure}[htbp] \centering \approxspace{ 4in} 
\includegraphics[angle=270]{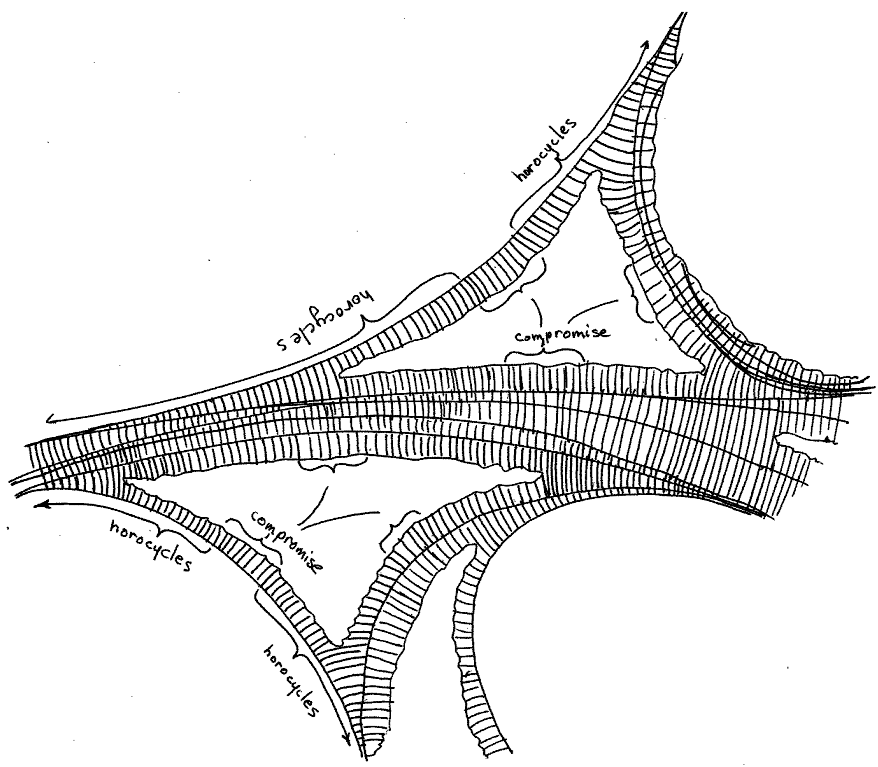}
\label{Construction of F(mu)} 
\caption{  
There is an orthogonal measured foliation $F ( \mu )$ defined in a small
neighborhood of any geodesic lamination $\mu$.
The transverse measure of $F ( \mu )$ agrees with arc length along the leaves
of $\mu$.  $F ( \mu )$ is canonical up to restrictions to smaller
neighborhoods and isotopies which are fixed along $\mu$.}
\end{figure}

A bounded total length of the sides of $S$ cut by $\mu$ remain bare:
this length consists of sides of simple closed leaves of $\mu$
together with a compact portion of each infinite side of $S$ cut by
$\mu$.  Near the closed leaves of $\mu$, let the leaves of $F(\mu)$
consist of the perpendicular geodesics.
Along the infinite sides, construct $F(\mu)$ by
compromising between the horocycles coming from one end and
the horocycles coming from the other, using an averaging procedure.
The foliation is so far constructed only in $N_{\epsilon} ( \mu ) - \mu$.
Its tangent line field is continuous at $\mu$, however, and in fact
its continuous extension is a Lipschitz line field.
Therefore, it is integrable, and it defines a foliation in all of
$N_{\epsilon} ( \mu )$.

A transverse measure for $F ( \mu )$ is defined by arc length along
leaves of $\mu$.  There is also a foliation $G ( \mu )$
whose leaves consist of orthogonal trajectories of $F ( \mu )$
and which contains $\mu$ as a closed invariant subset.  In most places,
the leaves of $G ( \mu )$ are geodesics, and transverse measure for 
$F ( \mu )$ agrees with arc length for $G ( \mu )$ --- everywhere except
in the compromise zone between two families of horocycles.

Observe that $F ( \mu )$ is essentially independent of the
choices: it is well-defined up to restrictions to a smaller neighborhood
followed by isotopies which are fixed on 
$\mu$.  Hence, the measure class of $F ( \mu )$
can serve as an invariant for the hyperbolic structure in a
neighborhood of $\mu$; when we want to study this dependence, we can
write it $F_{g} ( \mu )$, where $g$ is the hyperbolic structure.

In \S\ref{Derivatives of lengths of laminations}
we will compute the derivative of the length function on lamination
space in terms of $F ( \mu )$. In
\S\ref{Cataclysm coordinates for Teichmuller space}
we will further analyze foliations such as $F ( \mu )$, which
are defined in a neighborhood of $\mu$ and are transverse to $\mu$:
they are parametrized in a nice way by a convex
subset of Euclidean space and they serve to parametrize Teichm\"uller space.

For the present, the task is the inverse construction:
suppose that we have a measured foliation $F$ defined in some neighborhood of 
$\mu$, with leaves transverse to the leaves of $\mu$.
Can we construct a hyperbolic structure $h$ in a neighborhood of $\mu$ such that
$F_{h} ( \mu )$ is equivalent to $F$? What additional information do we need?

There is some additional information which can be easily read from a hyperbolic 
structure in a neighborhood of $\mu$.  In a neighborhood of any ideal vertices 
$q$ of the complement of $\mu$, the foliation is composed of horocycles.
In this neighborhood, there is a {\it sharpness}
function $f_{h,q}$, whose value at a point $x$ is
defined as the log of the length of the horocycle through 
$x$ about the given ideal vertex.  The absolute value of the 1-form 
$df_{h,q}$ is the transverse measure for $F ( \mu )$ in the given spike.
Given the measure class of $F ( \mu )$, this condition determines 
$f_{h,q}$ up to a constant.

\begin{proposition}[Neighborhood of a lamination]\label{Neighborhood of a lamination}
If $F$ is any measured foliation defined in $N_{\epsilon} ( \mu )$
with leaves transverse to the leaves of $\mu$ and if $F$ is standard in a
neighborhood of any cusps of $S$ (this means that its leaves are circles,
and there is infinite transverse distance to the cusp), then there is a
hyperbolic structure $h$ of finite area on the neighborhood such that
$F_{h} ( \mu ) = F$.  Furthemore, the cusps of $S$ are still cusps.

The set of all hyperbolic structures $h '$ of finite area on the neighborhood
having $F_{h}' ( \mu ) = F_{h} ( \mu )$ are classified
(up to restriction to smaller neighborhoods) by the constant values
$f_{h,q} -  f_{{h} ' , q }$, which can be chosen arbitrarily.
\end{proposition}

\proof %{Proof of \shortlabel Neighborhood of a lamination:4.3:}
Let $F$ be a measured foliation defined in a neighborhood $N$ of $\mu$
and satisfying the hypothesis of the proposition.  For each ideal vertex $q$
of the complement of $\mu$, the restriction of $F$ to a sufficiently small
neighborhood is equivalent to the foliation by horocycles of the region
between two parallel lines in the hyperbolic plane.
This is an elementary consequence of the fact that each
infinite end of a leaf of $\mu$ has infinite transverse measure.  Therefore,
we can choose a function $f_{q}$ in this neighborhood so that $f_{q}$
goes to $+\infty$ near $q$ and $| d f_{q} |$ agrees with the transverse
measure of $F$ .

There is now a unique hyperbolic structure $h$ defined in a neighborhood of $q$
such that $f_{h,q} = f_{q}$.  These hyperbolic structures extend in a unique
way to give a hyperbolic structure on $N - \mu$ such that along the edge of any
component of $N - \mu$, arc length agrees with transverse measure of 
$F$.

We will extend the hyperbolic structure which is currently
defined in the complement of $\mu$ to all of $N$ by
describing its developing map as an infinite product in the
group of isometries of the hyperbolic plane.

Pick a base point $x_{0}$ in $\hy^{2}$, and a tangent vector $V_{0}$ at 
$x_{0}$.  For any point $y \element N - \mu$, there is an isometry of
a neighborhood of $y$ to the hyperbolic plane taking $y$ to $x_{0}$
and a tangent vector of $F$ to $V_{0}$.

%\begin{figure}[htbp] \centering \approxspace{2in} 
%\includegraphics{dummy.eps}
%\label{developing condition} 
%\caption{The values of
%the developing map $D$ for a leaf of
%$F(\mu)$ at the two endpoints of an interval of the intersection with
%the complement of $\mu$ are related by a certain parabolic transformation}
%\end{figure}

Now we will construct a developing map $D$ along the $F$-leaf $m$ through
$y$.  This leaf intersects the complement of $\mu$ in a union of open intervals
$(a_{i} ,  b_{i} )$ .  Whatever $D$ is, it must satisfy the condition that
$$ D (b_i ) =  P_i D ( a_i )
$$
where $P_{i}$ is a parabolic transformation fixing one of the two endpoints
of the line perpendicular to $V_{0}$, depending on the direction to the
particular ideal vertex, and moving $x_{0}$ a distance
$\exp ( - f_{h,q} ( a_{i} ) )$ along a horocycle.
\begin{figure}[htbp] \centering \approxspace{ 5in} 
\includegraphics[angle=90]{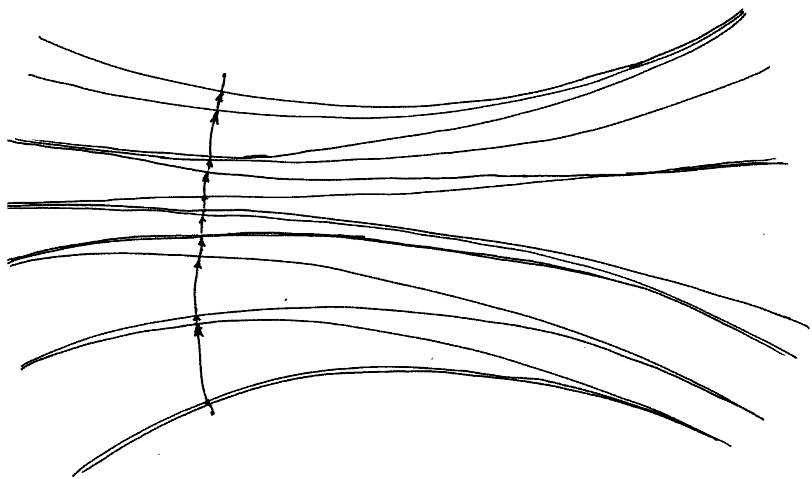}
\label{developing the neighborhood of a lamination} 
\caption{  
The developing map for a hyperbolic structure in a neighborhood
of a lamination can be expressed as an infinite composition of an
linearly ordered set of isometries,
with linear order depending on the topology of the lamination.  }
\end{figure}
If $y '$ is another point on $m$, we may define $D ( y ' )$
to be the product of the transformations $P_{i}$
associated with the intervals of $(y, y ' ) - \mu$, but first we must make
sense of the infinite product.  The $P_{i}$ form a linearly ordered set
according to the order of the intervals $(a_{i} , b_{i} )$,
but the linear order is generally
not the same as the linear order of a sequence of integers.
Order of multiplication is important in a noncommutative group, so we
must not disturb the given linear order.

The definition for convergence of a nonabelian linearly ordered countable
product is clear: the infinite product converges if
products over finite subsets converge to a limit as the finite set grows.

For any compact subset of the leaf $m$, there is a lower bound
to the distance measured within a particular spike
between any two of the components of intersection of 
$m$ with the spike.  The lengths of these horocycles therefore have a
sequence of lengths estimated from above by a geometric progression
(when they are reordered according to the value of 
$f_{h,q}$).
There is only a finite set of ideal vertices, so the series of
lengths of horocycles converges absolutely.

Choose any smooth left-invariant metric on the group of isometries
of the hyperbolic plane and write the distance of a group element $g$
from the identity as the norm $\|g\|$.  We can conclude that the series 
$K= \sum_i \| P_{i} \|$ converges absolutely.  By the triangle inequality,
the $K$ is an upper bound to the norm of any finite subproduct.

Consider the effect of inserting a small group element $e$ in the middle of
a finite product.  What is the distance between the two results?
We need to compare a product $A B$ with a product $A e B$.  We write
$$ A e B = A B ( e )(e \inverse B \inverse e B )
$$
Since the commutator $e \inverse B \inverse e B$ is the identity when
either $e$ or $B$ is the identity, it follows by an elementary
calculus argument that when $e$ and $B$ are restricted to the compact set
of elements with norm less than $K$, the norm of their commutator
is less than some constant times $ \| e \|  \| B \|$ .
Therefore, the distance from $A B$ to $A e B$ is less than some
constant times $\|e\|$.  The finite products satisfy the Cauchy criterion,
so they converge.

This definition of the developing map restricted to the leaf
extends in an obvious way for points on either side of the leaf,
to give a developing map for a neighborhood of the leaf $m$.

The finite products in this discussion have a geometric interpretation.
If the lamination $\mu$ is simplified in a neighborhood of $m$
by removing all but a finite number of the gaps between leaves of $\mu$,
thereby identifying the leaves of $\mu$ to a finite set of leaves,
a finite local approximation to $\mu$ is obtained.
Each finite product defines a developing map for a hyperbolic neighborhood of 
$m$ for some finite approximation to $\mu$.  Comparing with these finite
products, it is clear that the local hyperbolic structure does not depend on
the choice of $y$.  Therefore, we have constructed a hyperbolic neighborhood of 
$\mu$ with the given $f_{h,q}$ .

It is easy to see that near a cusp, the hyperbolic structure is
standard, on account of the condition on $F$.

To complete the proof of the proposition, we need to show that hyperbolic
neighborhoods are determined by $F_{h} ( \mu )$ and $f_{h,q}$.
Suppose that we have two such hyperbolic neighborhoods.  By hypothesis,
there is a map $H$ between them, preserving $\mu$ and the foliation $F$.
This map can be isotoped to be an isometry in the complement of $\mu$,
using the equality of the functions $f_{h,q}$.  It is essential that the
lamination $\mu$ have 2-dimensional Lebesgue measure 0 in both hyperbolic
neighborhoods.  The argument which works for the case of a closed surface
carries over to the present case:
the area of the neighborhood can be computed using the Gauss-Bonnet theorem,
from its Euler characteristic together with the total geodesic curvature
of the boundary of the neighborhood.
The total area of the complement of the lamination in this neighborhood
can also be computed from the Gauss-Bonnet theorem,
and it is the same as the area of the entire neighborhood.

It now follows that the arc length of a leaf of $F ( \mu )$
is the sum of the lengths of intersection with the complement of 
$\mu$, so arc length along leaves of $F ( \mu )$ agrees for the two
hyperbolic structures.  This implies the developing maps are conjugate by
an isometry, hence the metrics are the same.
\endproof %{\ref {Neighborhood of a lamination}}

\begin{corollary}[stretch defined]\label{stretch defined}
For any complete hyperbolic structure $g$ of finite area on a surface $S$,
for any complete geodesic lamination $\lambda$ not all of whose leaves go
to a cusp at both ends there is a new hyperbolic structure
$$ \Stretch( g , \lambda , t )
$$
depending analytically on $t > 0$ such that
\begin{itemize}
\item[(a)] the identity map is Lipschitz with Lipschitz constant
$$ L ( g , \Stretch ( g , \lambda , t ))  = b \exp ( t ) ,
$$
and\\
\item[(b)] the identity map exactly expands arc length of $\lambda$ by the constant
factor $\exp ( t )$.
\end{itemize}
\end{corollary}

\proof %{Proof of \shortlabel stretch defined:4.6:}
A new measured foliation is constructed in all of $S$
simply by multiplying the old measure by a constant.
Use proposition \ref {Triangle expands}
to define a new hyperbolic structure in each ideal triangle,
and use these structures to determine the new sharpness functions:
they are $\exp ( - t ) f_{g,q}$.  Apply the previous proposition to extend
the hyperbolic structure over a neighborhood of $\lambda$.

To see that the dependence on $t$ is analytic, we consider the holonomy of 
$\Stretch ( g , \lambda , t )$ around any element of the fundamental group
of the $S$.  Represent the element by a sequence of arcs along leaves of 
$F ( \lambda )$ interlaced with arcs along leaves of $\lambda$.
The net holonomy is expressed as an infinite product, by the preceding proof.
All the terms in the product, even those representing
the finite number of segments of the loop along $\lambda$, have the form
$\exp (  \exp ( \pm t )  V_{i)} $ where the outside $\exp$ represents
exponentiation in the Lie group, and $V_{i}$ is an infinitesimal isometry.
The series $\sum \| V_{i} \|$ converges absolutely, so the product converges
not only for all $t \element \reals$ but for all $t \element \complexes$.
(Note that the
product can still be given an interpretation in this case, since
the complexification of the group of isometries of $\hy^{2}$ is
the group of isometries of $\hy^{3}$: one gets developing maps of the
surface into $hy^3$ for complex $t$.)

This is a product of holomorphic functions of $t$,
so the product converges to a holomorphic function.
\endproof %{\ref {stretch defined}}

Note that most hyperbolic structures defined in a neighborhood of
the lamination $\lambda$ do not extend to complete hyperbolic surfaces:
in general, there is likely to be non-trivial holonomy going around the
boundary of an ideal triangle

\section {Convex models for measured lamination space} \label {Convex models for measured lamination space}

There are a number of analogies between the space $\PL ( S )$
of compactly supported projective laminations (nonzero measured
laminations up to positive real multiples) and the boundary of a convex set.
This is no accident.  Consider any projective lamination, represented by a
measured lamination $\mu$.  The logarithmic derivative of the length of 
$\mu$ is an element of the cotangent space of Teichm\"uller space $\teich(S)$
which depends only on the projective class of $\mu$.

\begin{theorem}[d log length is convex]\label{d log length is convex}
For any hyperbolic structure $g$ on a surface, the function
$$d \log\length : \PL(S) \arrow T_g^*(\teich(S))
$$
embeds the space of projective laminations on 
$S$ as a convex sphere containing the origin in
the cotangent space to Teichm\"uller space at $g$.
\end{theorem}

\proof %{Proof of \shortlabel d log length is convex:5.1:}
It is clear that $d\log\length$ defines a continuous map of
the sphere of projective laminations into the cotangent space of
Teichm\"uller space.

Let $\lambda$ be any measured lamination on $S$.  By adding a finite number of
leaves, $\lambda$ can be extended to a geodesic lamination $\mu$ whose
complement $S - \mu$ consists of ideal triangles and with no invariant
measure whose support is larger than the support of $\lambda$.
Consider the tangent vector $X$ to the curve $\Stretch ( g , \mu , t )$
in Teichm\"uller space.  The derivative of the log of the
length of $\lambda$ in the direction $X$ is 1, while the derivative of the
log of the length of any other measured lamination does not exceed 1.
In other words, $X$ defines a linear function on the cotangent space which
has maximal value on the image of $\PL ( S )$ at $d \log \length ( \lambda )$,
proving that the image of $\lambda$ is an extreme point of the convex hull
of the image of $\PL ( S )$.  (An extreme point of a convex set is a point for
which there is a linear function which is not constant on the set which attains 
its maximum at the given point.) This linear function also shows that the
image of $\lambda$ is distinct from the image of any other projective
lamination except possibly other laminations represented by some different
invariant measure $\lambda '$ on the same geodesic lamination.

To show that $\PL ( S )$ is really embedded, consider what happens to length of 
$\lambda$ and $\lambda '$ under small twists along
simple closed curves nearly tangent to the underlying geodesic lamination.
The lengths change approximately in proportion to the intersection numbers
with the simple closed curve.  It is easy to find nearly tangent curves
with different intersection numbers, so $d\log\length$ is indeed an embedding.

The set of extreme points of any convex set is a sphere of some dimension.
The dimension of the Teichm\"uller space of $S$ is one more than the dimension
of $\PL ( S )$, so the only possibility
is that the set of extreme points of the convex hull of its image
must coincide with the image, since a sphere is not homeomorphic to
any subset of a lower dimensional sphere and not homeomorphic to a
proper subset of itself.

Since the convex hull has an interior,
there must be at least one line through the origin which intersects
the image of $\PL ( S )$ in at least two points.
For each of these points, there is a linear function which
attains its positive maximum value there.
It follows that 0 must separate the two points, so 0
is a convex combination of them.
\endproof %{\ref {d log length is convex}}

Two qualitative features of the convex embedding should be noted.
The first feature comes from any geodesic lamination $\mu$ such that 
$\mu$ admits more than one projective class of invariant measures.
These exist for any surface more complicated than a punctured torus;
the easiest examples are collections of simple closed curves,
but there are even examples for surfaces of genus 2 or more where
the lamination itself is minimal.
The examples are called non-uniquely ergodic laminations; in general, the
set of projective classes of invariant measures is a simplex in projective
lamination space.
The length for a measure is a linear function of the measure, so each such
simplex appears as a flat simplex on the image of the projective
lamination space.
\begin{figure}[htbp] \centering \approxspace{4in} 
\includegraphics{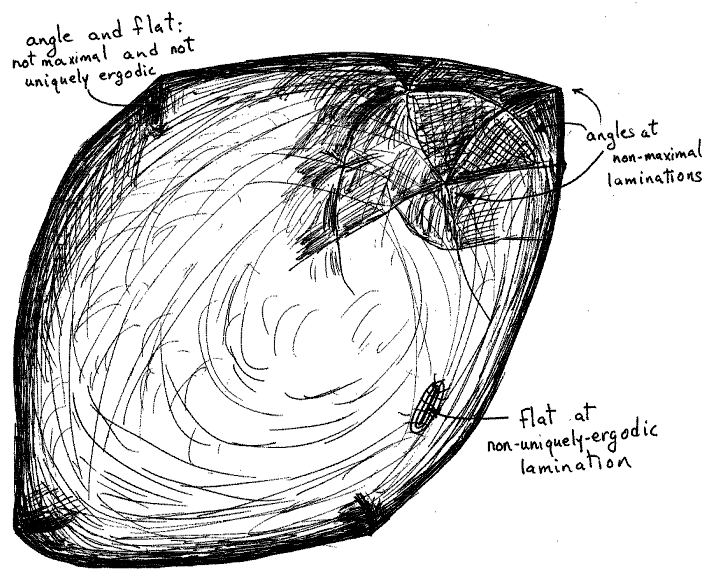}
\label{ convex picture  of projective laminations } 
\caption{
In the convex image of projective lamination space in the cotangent
space to Teichm\"uller space, the image of any nonuniquely ergodic
lamination is on a flat, and the image of any non-maximal compactly supported
lamination is on an angle.
This three-dimensional picture is not meant literally, since the
first oriented examples which illustrate both features have four-dimensional
Teichm\"uller spaces: the twice-punctured torus or the quintuply-punctured
sphere.  }
\end{figure}

The second feature comes from projective laminations which are not maximal
as geodesic laminations.
For any extension to a maximal lamination (with complementary regions
consisting of ideal triangles), a stretching vector is defined,
giving a linear function attaining its maximal value at the given
point.  Provided these vectors are truly distinct, this means that there is
an angle in the convex sphere at the given point.
This phenomenon already occurs for the punctured torus, where there is
an angle at any projective lamination which is a simple closed curve.

The convexity can also be expressed in terms of the space $\ML$ of
measured laminations.  The 1-jet of the length of a measured lamination is
described by the pair $( \length , d  \length )$.
This maps $\ML$ into a vector space of one higher dimension, and
its image is the boundary of a convex cone.
\begin{figure}[htbp] \centering \approxspace{6.7in} 
\includegraphics{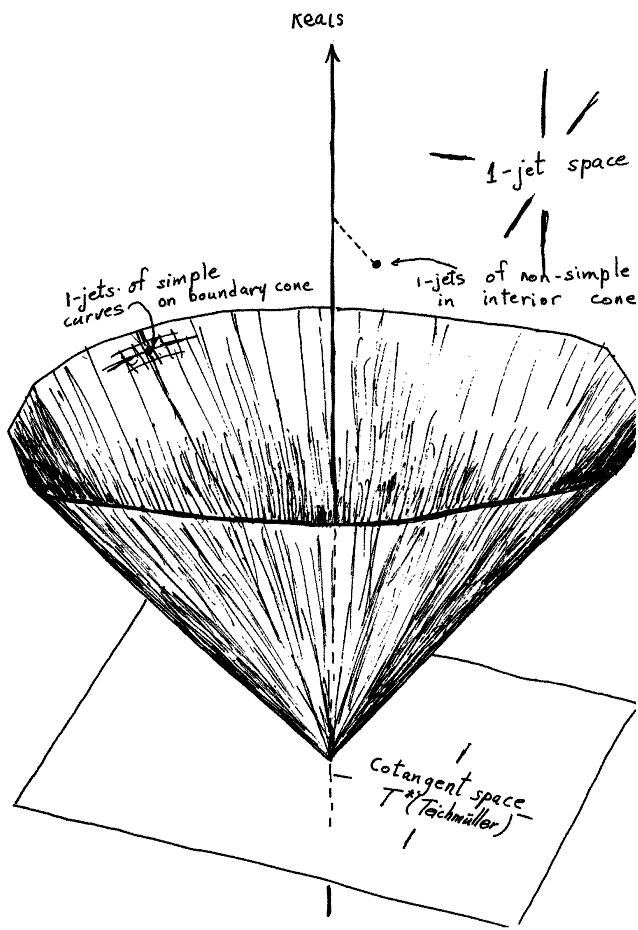}
\label{1-jet space} 
\caption{
The function which takes a measured lamination to the 1-jet of its
length at a point in Teichm\"uller space
embeds it as the boundary of a convex
cone in the space of 1-jets to Teichm\"uller space.
}
\end{figure}

A convex subset of the cotangent space determines a dual convex
subset of the tangent space.
This also determines a nonsymmetric norm, or nonsymmetric Finsler metric,
which is the infinitesimal version of $K$:
$$ \|Y\|_{\length} = \sup_\lambda
\left (  \frac{Y ( \length (\lambda) ) }{ \length ( \lambda ) } \right )
$$

In this dual convex set, there are flats which correspond to angles
in the original set, and angles which correspond to flats in the
original.

Associated with a measured geodesic lamination $\lambda$
is a certain vector field $E_{\lambda}$ on Teichm\"uller space,
its earthquake vector field.
There is a simple formula for the effect of an earthquake on the length
of a geodesic lamination, which has been used by Steve Kerckhoff
and Scott Wolpert (\cite{Kerckhoff}, \cite {Wolpert6}):
it is the integral over the set of intersections
of the cosine of the angle of intersection, with respect to the product
measure.
What is cogent here is a consequence, that for any two measured geodesic
laminations $\lambda$ and $\mu$, there is an antisymmetry
$$ E_\lambda ( \length ( \mu ) ) = - E_\mu ( \length ( \lambda ) )
$$

This formula should suggest
a nonsingular bilinear pairing on the tangent space of a manifold,
which determines
an identification of the tangent space with the cotangent space.
In the present circumstance,
we know maps of measured laminations both to the tangent
space (through earthquakes) and to the
cotangent space (through d length).
Given that at least one of the maps is a homeomorphism,
the formula implies that the other map is also a homeomorphism,
and that the induced linear structures on measured lamination
space agree, so we have an identification of the two vector spaces via
the space of measured laminations.
(The linear structure induced on measured lamination
space varies from point to point, though).

In fact, Scott Wolpert \cite{Wolpert1}
has shown that this antisymmetry can be interpreted
as the K\"ahler form for the Weil-Peterson metric, applied to two
the vectors $E_{\lambda}$ and $E_{\mu}$.

\begin{theorem}[Earthquakes are convex]\label{Earthquakes are convex}
The unit earthquake vector fields 
$$ \frac 1{ \length ( \lambda ) }  E_\lambda
$$
define an embedding of $\PL ( S )$ as a convex subset of the tangent space at
any point $g$ in Teichm\"uller space.
\end{theorem}

\proof %{Proof of \shortlabel Earthquakes are convex:5.4:}
If $\lambda$
is any measured geodesic lamination, then as we have seen there is some
tangent vector $X \in T(\teich(S))$ such that 
$X ( \log ( \length ( \lambda ) )$ takes its maximum value at $\lambda$.
By the earthquake theorem, there is some measured geodesic lamination $\mu$
such that $X= E_{\mu}$.  Consider $- d \log ( \length ( \mu ) )$
as a linear function on the tangent space to Teichm\"uller space.
From the antisymmetry formula, this linear function takes its maximum
value on the image of $\PL$ at $E_{\lambda}$.
\endproof %{\ref {Earthquakes are convex}}

In light of this fact, there is another nonsymmetric Finsler metric on
Teichm\"uller space, the earthquake norm, whose unit sphere is the set of
unit speed earthquakes.
Earthquakes are {\it not} geodesics for this norm ---
the earthquake vector fields are not (in some intuitive sense) ``parallel''.
The distance between two elements of Teichm\"uller space in the
earthquake norm is the infimum of the total magnitude of a sequence
of earthquakes transforming one to the other.
To make such a transformation,
it is highly inefficient to persist
with earthquakes along any particular fault for very long.
This can be
seen clearly for the punctured torus, whose Teichm\"uller space
equipped with the Teichm\"uller metric is isometric to the hyperbolic plane.
Earthquakes along simple curves are approximately
horocycles in this case; in particular, both ends limit to the
same point on $S_\infty^1$.

Perhaps there is some global minimization problem which the earthquake
norm measures, but I don't know any better description.
In contrast to the length norm,
I don't know any global interpretation for distances with respect
to the earthquake norm.

\section {The tangent space of measured lamination space} \label {The tangent space of measured lamination space}

We have now developed a good picture for the behavior
of lengths of simple closed curves and laminations
as a hyperbolic structure varies infinitesimally in Teichm\"uller space,
in terms of the length norm, a nonsymmetric Finsler metric.

The global problem still remains, however.
Even though we have constructed some long-range stretch paths which
are geodesics for the the length norm, there are definitely
gaps between the stretch paths.  We cannot pass between an arbitrary pair of
hyperbolic surfaces with a stretch path.
What we will show is that any two hyperbolic structures
can be connected by geodesics for the Finsler metric which are composed of
a finite sequence of the $\Stretch$ paths, bounded in number.

The phenomenon that geodesics do not have unique extensions, and in fact
can turn through limited angles if they want, is related to the fact
that the unit ball for the Finsler metric is not strictly convex.
A simple example for comparison is the $L^{1}$ or
Manhattan metric in the plane, given by the norm
$$ \| ( x ,  y ) \|_M  = |x| + |y|
$$
on the tangent space to the plane.
Any curve whose tangent vector is always in a single quadrant is
a geodesic.
Any two points can be connected by a geodesic which first goes along
the $x$-axis, then along the $y$-axis.  The San Francisco metric
$$ \| ( x , y ) \|_{SF} = dh ( x , y ) + |x| + |y|
$$
is a simple variation which is nonsymmetric, where $h$
is the height function.  It is assumed that $h$
has first derivative of norm less than 1 everywhere. Unfortunately,
I do not have the analytic expression for the height function $h$
of San Francisco.  Note that the $\| \|_{SF}$-geodesics are also
$\| \|_M$-geodesics, since the integral of $dh$ is
independent of path.
There are many variants, depending on preferences for locomotion.

In this section, we shall analyze the tangent space to measured lamination
space, insofar as this makes sense.
After we have developed the appropriate theory, we will be able,
in section
\ref{Derivatives of lengths of laminations},
to compute the derivative of the length functions within
measured lamination space for a fixed hyperbolic structure.
This computation will connect the topologically defined
tangent space of measured lamination space to the 
geometrically defined tangent space (or tangent cone)
of a convex body.

In addition, just as the behavior of lengths of laminations
infinitesimally in Teichm\"uller space gives a good global picture
for measured lamination space, so we shall see that the behavior of lengths of
laminations infinitesimally in measured lamination space will
give good global coordinates for Teichm\"uller space.

$\ML ( S )$ has a natural piecewise linear structure.
The nicest way to see this structure seems to be in terms of train tracks.
Any maximal recurrent and transversely recurrent train track defines a closed
cell of maximal dimension in the measured lamination space.
(Recurrent means that there is a closed loop through every branch of the
train track.  Transversely recurrent means there is a closed loop transverse to 
the train track which never ``doubles back'' passing through any point.
Transversely recurrent train tracks are precisely those which can be
embedded in some hyperbolic surface so that the branches all have
curvature less than $\epsilon$ for any $\epsilon$.
Transverse recurrence is a technical convenience which makes for a smoother
development of the theory.
A maximal train track is one whose complementary regions are all triangles
and punctured monogons.)
The coordinate changes when two of these cells intersect are always piecewise
linear.

Measured lamination space, like any other piecewise linear space,
has a tangent space associated with any point, made up of tangent
vectors to simplices at that point.  The tangent space to a piecewise
linear space does not have a natural linear structure, but
the operation of multiplication by scalars is well-defined.
In other words, the structure is that of a cone.
Any piecewise linear map induces a piecewise linear map from the tangent space
at any point to the tangent space at its image.

Unfortunately, the tangent spaces at the points in a piecewise linear
space do not piece together to form a tangent bundle.
There is no way to define a topology on the union of the tangent spaces
to form a bundle on which the induced transformations are continuous.
For instance, a piecewise-linear function $\reals\arrow\reals$
has a derivative which is a step function, and resists all attempts to
make it continuous.

There is considerably more structure on measured lamination space than
just this piecewise linearity, however.
Non-linearity can occur only along very special cells,
which are defined by laminations carried on recurrent and transversely recurrent
train tracks which are not maximal.
At any point not on one of the countable set of special cells
there is a well-defined linear structure on the tangent space:

\begin{proposition}[Generic tangent space is linear]\label{Generic tangent space is linear}
If $S$ is a complete hyperbolic surface of finite area and if $\lambda$
 is any complete measured lamination of compact support
(that is, its complementary regions
are all ideal triangles or punctured monogons) then the space of
compactly supported measured lamination has a well-defined linear structure at 
$\lambda$, inherited from any train track which carries $\lambda$.
\end{proposition}

\proof %{Proof of \shortlabel Generic tangent space is linear:6.1:}
If $\tau_{1}$ and $\tau_{2}$ are any two train tracks which carry $\lambda$,
there is a third train track $\sigma$ which carries $\lambda$ and is carried by 
$\tau_{1}$ and $\tau_{2}$.  In fact, any $\tau_{\epsilon} ( \lambda ) $
with $\epsilon $ sufficiently small may be taken as $\sigma$.
Each of the three train tracks give coordinate systems for an open
neighborhood of $\lambda$ (since $\lambda$ is a complete lamination),
and the maps from $\sigma$ to the $\tau_{i}$ are linear.  It follows that
the linear structure of the tangent space to $\lambda$
inherited from any train track carrying it is the same linear
structure it inherits from any other.
\endproof %{\ref {Generic tangent space is linear}}

At a lamination $\lambda$ which is not complete, the tangent space to measured
lamination space does not in general have a linear structure.
Still, fragments of a linear structure remain, and it is possible to describe
the tangent space as a union of these fragments, where each linear fragment
is labeled by a certain kind of lamination which contains 
$\lambda$.

To begin the analysis, note that any measured lamination $\lambda$ has a
canonically associated set $C ( \lambda ) $ of topological laminations,
consisting of the intersection of the closures
in the Hausdorff topology
\footnote{The Hausdorff distance between two closed subsets
of a metric space $X$ is the maximum distance from a point in either set to
the nearest point of the other set; this defines a topology on the
set of closed subsets independent of metric if $X$ is compact,
and in general depends at most on the Lipschitz-class
of the metric.}
of its neighborhoods in the
measure topology.  The condition that $\lambda$ be complete is precisely the
condition that $C ( \lambda ) $ consist of only $\lambda$ itself.

What does $C ( \lambda )$ look like in general?
A necessary condition is given in terms of the following definition:
A lamination $\mu$ is {\it chain recurrent} if for any
$\epsilon$ and for any point $x \element \mu$ there is a closed 
$\epsilon$-trajectory of $\mu$ through $x$, that is, a closed unit speed path
in the surface such that for any interval of length $1$
on the path there is an interval of length $1$ on some leaf of $\mu$
such that the two paths remain within $\epsilon$ of each other in the 
$C^{1}$ sense.

Any measured lamination is recurrent, and any recurrent lamination
is chain recurrent.  The relevance of the condition of chain recurrence is that,
unlike the condition of recurrence, it
is a closed condition in the Hausdorff topology on laminations:

\begin{proposition}[chain recurrence is closed]
\label{chain recurrence is closed}
Any compactly-supported lamination which is a Hausdorff limit of chain
recurrent laminations is chain recurrent.
\end{proposition}

\proof %{Proof of \shortlabel chain recurrence is closed:6.2:}
Given $\epsilon$, there is a neighborhood of any compactly supported lamination 
$\mu$ in the Hausdorff topology such that an $\epsilon/2$-trajectory for any
lamination in that neighborhood is an $\epsilon$-trajectory for $\mu$.
\endproof %{\ref {chain recurrence is closed}}

\begin{proposition}[Hausdorff closure of neighborhoods]\label{Hausdorff closure of neighborhoods}
$C ( \lambda )$ consists of all chain recurrent laminations which contain 
$\lambda$.
\end{proposition}

\proof %{Proof of \shortlabel Hausdorff closure of neighborhoods:6.3:}
Given any chain recurrent lamination $\mu$ containing $\lambda$,
consider any $\epsilon$-train track approximation $\tau$ to $\mu$.
Clearly, $\tau$ is recurrent, since a closed near-trajectory of 
$\mu$ gives a closed trajectory on $\tau$.
Therefore, it carries some measured lamination $\nu$ which covers every branch
of $\tau$.  A convex combination of the weights for $\lambda$ and $\mu$
defines a measured lamination $\nu$ which is arbitrarily near to 
$\lambda$ in the measure topology, and is carried on $\tau$
with positive weights on each branch.

If this construction is repeated for a sequence of finer and finer train
track approximations to $\mu$, the limit of the associated $\nu$'s is
$\lambda$ in the measure topology but $\mu$ in the Hausdorff topology.
\endproof %{\ref {Hausdorff closure of neighborhoods}}

A tangent vector to measured lamination space at $\lambda$ needs to describe
more data than a Hausdorff limit lamination:
it must also describe limiting behavior of measures.  The laminations in 
$C ( \lambda )$ all contain $\lambda$, so they
are determined by their behavior in the complement of $\lambda$.
Define $S_{\lambda}$ to be the hyperbolic surface obtained from $S - \lambda$
by metric completion.

$S_{\lambda}$ is a hyperbolic surface with geodesic boundary which
may have components isometric to a line, linked together in chains
like an ideal polygon.  Let $\far_{\lambda} $ be the space of measured
laminations on $S_{\lambda}$, where leaves are not allowed to intersect
the boundary, but they may limit on ideal vertices of $S_{\lambda}$
or spiral around a closed boundary component in either direction.

For $\lambda '$ in a small neighborhood $U$ of $\lambda$ in $\ML ( S )$, leaves
of $\lambda ' $ can cross leaves of $\lambda$ only at small angles.

A map
$$ \cut_\lambda : U\arrow \far_\lambda
$$
is defined by ignoring all leaves of $\lambda '$
inside $N_{\epsilon} (\lambda)$,
pushing all leaves of $\lambda '$ which intersect an infinite $\lambda$-leaf to
the ``nearby'' ideal vertex of $S_{\lambda}$, and spinning any
$\lambda'$-leaf which
cross a closed $\lambda$-leaf about the closed leaf.  Some leaves of $\lambda'$
may coalesce in the process; still, the transverse measure for
$\lambda'$ pushes forward to give a
transverse measure for $\cut_\lambda(\lambda')$.

Laminations in the space $\far_{\lambda}$ are made up of two essentially
different parts.  First, there are the remote components which do not
intersect a small neighborhood of $\lambda$,
{\it i.e.}, the compactly supported laminations in $S_{\lambda}$.
Second, there are
a finite number of leaves with both ends tending toward 
$\lambda$.  (There can only be a finite number of leaves of $\cut(\lambda')$
with both ends tending
toward $\lambda$, because such leaves are determined by their limiting behavior
at the two ends, and their homotopy class.  There is a bound to the number
of homotopy classes represented by non-intersecting arcs on any
finitely-connected surface.  There can be no leaves of $\cut(\lambda')$
which tend at only one end toward $\lambda$, because each
component of $\lambda'$ is minimal, so that a leaf of $\lambda'$
which comes near $\lambda$ comes near infinitely many times.)

The map $\cut_{\lambda}$ also makes sense for geodesic laminations
(without transverse measures) which intersect $\lambda$ only at
small angles; in this case, of course, they are mapped to laminations
on $S_{\lambda}$ without transverse measures.

\begin{figure}[htbp] \centering \approxspace{ 4in} 
\includegraphics{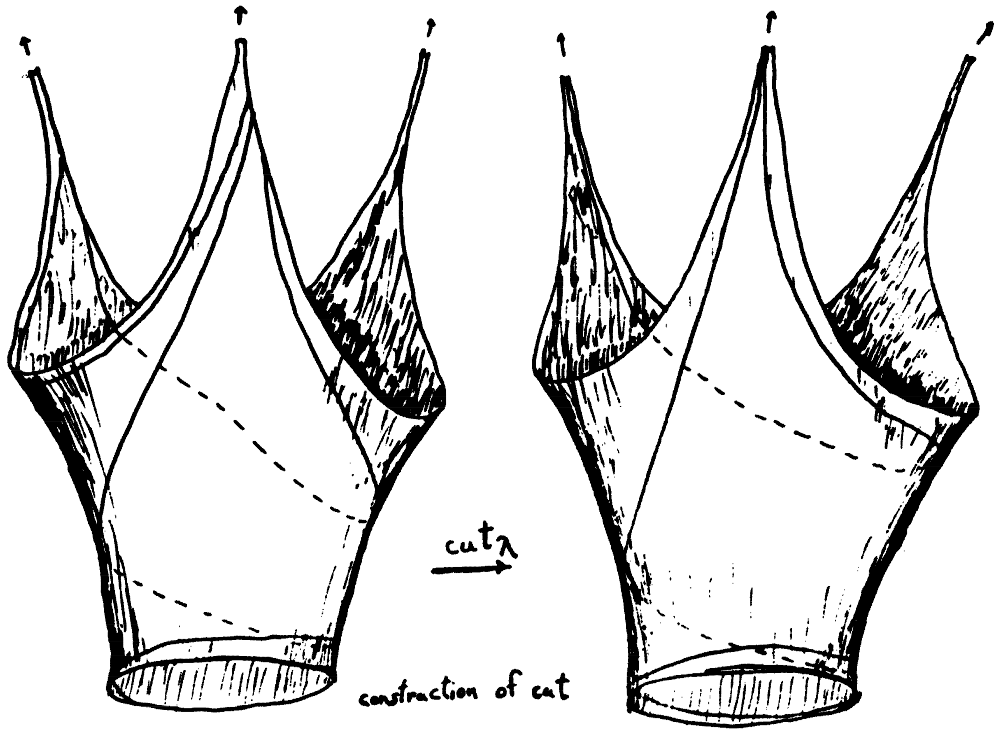}
\label{ construction of  cut} 
\caption{ 
The map $\cut_{\lambda}$ assigns to a lamination which does not
cross $\lambda$ at a large angle a new lamination in the surface cut
open along $\lambda$.}
\end{figure}

The definition of the map $\cut_{\lambda}$ can even be extended to all
laminations which do not intersect $\lambda$ at right angles anywhere,
but the definition then depends grossly on the particular hyperbolic surface
$S$.  For any two hyperbolic structures, the definitions agree
for sufficiently small neighborhoods $U$.

This construction carries over to train track coordinates.
Consider any chain recurrent lamination $\mu$ containing $\lambda$.
Let $\sigma$ be a fine train track approximation of $\mu$, fine enough
that a regular neighborhood of $\sigma$ has the stable topological type
of small neighborhoods of $\mu$, and also the image of $\lambda$ on $\sigma$
is a sub-train track $\rho$ whose regular neighborhoods
have the stable topological type of a
small neighborhood of $\lambda$.

There is a quotient space of the space of invariant weights on $\sigma$,
obtained by forgetting the weights on $\rho$.
If $\mu$ is represented by a set of weights $m$ on $\sigma$,
then the image $\cut_\lambda(\mu)$ depends only on the image of $m$ modulo
weights on $\rho$.
Let $W_{\sigma} ( \mu )$ be the subset of the quotient space (invariant
weights on $\sigma$ modulo those which are zero outside $\rho$)
parametrizing
laminations which have the topological type of $\mu - \lambda$.
Let $X_{\sigma} ( \mu )$ be the set of weights on $\sigma$ which map to
$W_{\sigma}  ( \mu )$.

\begin{proposition}[Linear fragments from chain recurrence]\label{Linear fragments from chain recurrence}
Let $S$ be any complete hyperbolic surface of finite area, and
$\lambda$ any measured lamination of compact support on $S$.
For any $\mu \element C ( \lambda )$ such that the lamination
$\cut ( \mu )$ admits a transverse invariant measure of full support,
there is a convex cone $V ( \mu )$ in the tangent space at $\lambda$,
which agrees with the cone of laminations parametrized by
$X_{\sigma} ( \mu )$
in any sufficiently fine train track approximation of $\mu$.

The tangent space at $\lambda$ is the union of the $V ( \mu ) $.
\end{proposition}

\proof %{Proof of \shortlabel Linear fragments from chain recurrence:6.5:}
First we check that the sets $X_{\sigma} ( \mu )$ are independent of the
choice of $\sigma$, up to the isomorphism inherited from the canonical
linear maps induced when one train track carries another, provided 
$\sigma$ is a sufficiently fine approximation to $\mu$.
Also, we check that the linear structure is independent of the choice of 
$\sigma$.

The linear structure is no problem: the set of $\epsilon$-train track
approximations $\tau_{\epsilon} ( \mu )$ to $\lambda$ are cofinal
in the set of all train tracks carrying $\mu$.
When $\epsilon_{1} < \epsilon_{2}$, then the cell determined by
$\tau_{{\epsilon}_{1}} ( \mu )$ maps linearly to the cell determined by
$\tau_{\epsilon_2} ( \mu )$.

To verify that the actual set $X_{\sigma}(\mu)$ is independent of choice,
consider any train track approximation $\sigma$ of $\mu$ which has
a regular neighborhood $N$ of the stable topological type.
We may suppose that the boundary of $N$ is formed from smooth curves,
meeting near the switches of $\sigma$ at $0$ angles, forming either
a smooth joint or a cusp as appropriate.  The lamination $\cut_{\lambda}(\nu)$
is defined by the slight adjustment of all leaves of $\nu$
which are not contained in the regular neighborhood of $\lambda$
itself so that they limit on $\lambda$ but do not intersect it.

For any sufficiently fine train track approximation, it is clear that the
set $X_{\sigma} ( \mu )$ parametrizes $V(\mu)$.
Therefore, $V ( \mu )$ is well-defined.

Now we must make sure that the $V ( \mu )$ cover the tangent space.
(Note that they do not form an open cover - many of the cones have a
relatively low dimension.)
What this means, in terms of the definition of the tangent space
for a piecewise linear space, is that for any piecewise linear
path $\lambda_{t}$ beginning at $\lambda$, there must be an interval
around zero which is contained in some $X_{\sigma} ( \mu )$.

The compactness of the space
of compactly supported laminations in the Hausdorff topology
guarantees that we can find at least one lamination 
$\mu$ which is a limit point of $\lambda_{t}$.
Clearly  $\mu$ is a chain recurrent lamination containing $\lambda$.
Let $\sigma$ be a fine train track approximation of $\mu$.
Then $\lambda_{t}$ is represented as a linear path, for sufficiently small
$t$.  The lamination $\cut_{\lambda} ( \lambda_{t}/t)$ is actually constant.
\endproof %{\ref {Linear fragments from chain recurrence}}
It is good to keep in mind the comparison with the
picture of a convex body in $\reals^{n}$.
Associated with any point on the boundary of a convex body is a tangent
cone, which is a union of rays through the point.  The tangent cone
may be thought of as the boundary of the convex body obtained as the
limit of a sequence of linear expansions of the original, fixing the point
in question.
It may alternatively be thought of as the boundary of the intersection
of all half-spaces containing the convex body with defining plane passing
through the point.  At a generic point, the tangent space is linear.
At special points, the tangent space is not linear, but still for a
generic tangent ray through the nongeneric point, there is a well-defined
tangent plane.
There is a whole hierarchy of linear structures of this type, obtained
by taking the tangent cone, then if this is not linear the tangent cones
of the tangent cone at the various points on the tangent cone, etc.

The elements of $C ( \lambda )$, in this comparison, correspond to
the linear pieces which occur anywhere in this hierarchy, not just at the
top level.

The theorem means that there is an infinitesimal version of $\cut_{\lambda}$,
$$ \Cut_\lambda : T_\lambda ( \ML ( S )\arrow \far_\lambda
$$
where 
$$ \Cut_\lambda ( \frac d{dt} ( \mu ( t ) ) ) =
\lim_{t\arrow 0 }  \frac 1t \cut_\lambda ( \mu ( t ) )
$$
The kernel is a subspace $\Near_{\lambda} \subset  T_{\lambda} ( \ML ( S ) )$
consisting of the tangent vectors to the space of laminations
which map to 0 under $\cut_{\lambda}$, that is, laminations
which are contained in a small neighborhood of $\lambda$ on $S$.

$\Near_{\lambda}$ has a linear structure;
it is the intersection of all the cones
$V_{\mu}$.  The images of $\cut_{\lambda}$ and $\Cut_{\lambda}$ are in general
not the entire space $\far_{\lambda}$.  Whenever $\lambda$ has a closed
geodesic, any spiralling to the corresponding pair of boundary components of the
cut
surface must be in opposite directions and the total measure of leaves tending
to the two components must be equal.  More generally, whenever a component of 
$\lambda$ is orientable, the total measure of leaves
tending to that component in the two directions must be equal.
Using train track coordinates, it is easy to verify that this condition
is necessary and sufficient for an element of 
$\far_{\lambda}$ to be in the image.

\section {Derivatives of lengths of laminations} \label {Derivatives of lengths of laminations}

The length of a lamination is a continuous function on the product
of Teichm\"uller space with measured lamination space.
For any particular lamination, length is a real analytic function on
Teichm\"uller space, but how nice is the dependency on the lamination
when a hyperbolic structure is fixed? We will show that length
has a first derivative which is linear on each of the linear fragments 
$V_{\mu}$ of the tangent space to measured lamination space,
and we will compute the first derivative by a simple construction
in terms of a the perpendicular foliation $F ( \mu )$ associated with $\mu$.
What is an appropriate definition of first derivative in the PL context?
The definition is obvious when dealing only with piecewise-linear functions,
but we need a definition with greater generality.  A neighborhood of the
origin in the tangent space to a piecewise-linear space at a point $p$
has a well-defined exponential map to a neighborhood of $p$,
up to restriction to smaller neighborhoods.
We can define a function $f$ on the piecewise-linear space to have a
{\it first derivative} at $p$ if there is a continuous function $g$
defined on the tangent space at $p$ which is linear on rays and
acts as a first-order approximation to $f$ in a neighborhood of $p$:
that is, the difference of $f$ with $g$ transported by the exponential map
is $o(r)$, where $r$ is distance from the $p$ in any piecewise-linear metric.

Let $\mu$ be any lamination on a hyperbolic surface $S$, and let $\epsilon$
be sufficiently small that the topological type of $N_{\epsilon} ( \mu )$ is
stable.  Recall that the perpendicular foliation
$F ( \mu )$ is well-defined, up to restriction to a smaller neighborhood.
For any measured lamination $\lambda$ whose support is contained in
$N_{\epsilon} ( \mu )$, there is an intersection number
$i ( \lambda , F ( \mu ) )$.

There is a convenient way to compute this interection number using train tracks.
The $\epsilon$ train track approximation to $\mu$ can be mapped into 
$N_{\epsilon} ( \mu )$ so that all the branches are transverse to 
$F ( \mu )$.  Each branch $b$ then has some intersection number $l ( b )$
with $F ( \mu )$.  We think of $l$ as an element of the vector space 
$B$ spanned by the branches.  The intersection numbers are determined only
when an explicit map is chosen.  Since the map is canonical only up to
isotopy, the vector $l$ is well-determined only up to the subspace 
$SR \subset  B$ spanned by the switch relation vectors $v_{1} + v_{2} - v_{3}$
where the $v_{i}$ are functions which assign $1$ to the branch $b_{i}$, 0
elsewhere, and $b_{1} , b_{2} $ on one side and $b_{3}$ on the other are
three branches incident to a switch.

Any measured lamination $\lambda$ supported in $N_{\epsilon} ( \mu )$
is determined by a set of weights on the branches, represented by
a vector $w ( \lambda ) \element B$.  The condition that 
$w ( \lambda )$ be invariant is precisely the condition that it lie in the
orthogonal complement $SR^\perp $ of the switch relation vectors,
with respect to the natural inner product on $B$.

The inner product on $B$ induces a pairing between $SR^\perp$
and $B / SR$.  The intersection number $i ( \lambda , F ( \mu ) )$
is precisely the value of this pairing on $w ( \lambda )$ and $l$.

\begin{theorem}[Derivative is intersection number]\label{Derivative is intersection number}
Let $S$ be any fixed complete hyperbolic surface of finite area, and 
$\lambda$ any measured lamination of compact support on $S$.
The length function for measured laminations on $S$ has a first derivative at 
$\lambda$ which is linear on every linear fragment $V_{\mu}$ of the tangent
space $T_{\lambda} ( \ML ( S ))$.  Explicitly, the first-order affine
approximation to lengths of laminations $\lambda '$ such that
$\cut_{\lambda} ( \lambda ' ) = \cut_{\lambda} ( \mu )$
(as laminations without invariant measures) is given by the formula
$$ \length ( \lambda' )  \approx  i ( \lambda' , F ( mu ) ).
$$
The first order approximation is uniform for $S$ in any compact set of
Teichm\"uller space.
\end{theorem}

Note that when we apply the statement to the special case that
$\lambda$ is the empty lamination $0$, the linear fragment $V_\mu$ consist of
laminations $\lambda'$ defined by transverse invariant measure for
underlying topological lamination of $\mu$. In this case, the statements
yields the exact formula $i ( \lambda ' , F ( \mu )) = \length(\lambda')$:
since $\length( t \lambda')$ is linear in $t$, it must agree with its
derivative at $0$.
This consequence is also immediate from geometric considerations,
since arc length matches transverse measure of $F$ along the leaves of $\mu$.

We shall prove in section
\ref {Cataclysm coordinates for Teichmuller space}
that the first derivative depends on the hyperbolic structure for $S$
as a real analytic function in Teichm\"uller space.  In 
\S\ref {The maximum stretch lamination is almost always a curve} we will derive a formula for the remainder term for the first order
approximation, in terms of the geometry of the countable set of laminations
where there are extra saddle connections, along a piecewise linear path in
measured lamination space.

\proof %{Proof of \shortlabel Derivative is intersection number:7.1:}
We can analyze the behavior of laminations near $\lambda$ in the
measure topology using a fixed cover of a neighborhood of $\lambda$ by a
finite set of train track approximations to all chain recurrent laminations
containing $\lambda$.  For these train track approximations, we can assume
that there is a fixed subtrack $\rho$ common to all of them which
is the image of $\lambda$ and has the stable topological type.
Let $\tau$ be one of these train tracks, and $\mu$ be a chain recurrent
lamination carried by $\tau$.  For simplicity, we assume that each switch of 
$\tau$ is generic, that is, it has one branch
attached in one direction and two attached in the other direction.
For any train track,
there is some small change (by squeezing together branches) which yields
a new train track satisfying the genericity condition and still carrying
all laminations which it previously carried.

The idea of the proof will be to show
that a lamination $\lambda '$ sufficiently near $\lambda$ and supported in a
neighborhood of $\mu$ has most of its leaves quite near to $\mu$ in the 
$C^{1}$ topology, enough so that its length is closely estimated
by transverse measure of $F ( \mu )$.

The error is measured by the ratio of the transverse measure of 
$F ( \mu ) $ to arc length along $\lambda '$, and it differs from 1 for two
reasons. 

First, the leaves of $F$ spread apart somewhat as they move away from
$\mu$ in certain places, the compromise zones.
This effect depends quadratically on distance from $\mu$, since the leaves of 
$F ( \mu )$ are orthogonal to $\lambda$.

Second, $\lambda '$ crosses the leaves of $F ( \mu )$ not quite orthogonally.
The angle is estimated by the distance from $\lambda$ in the $C^{1}$
metric, and contribution to the error is to multiply the ratio by
the cosine of the angle.  Again, this depends quadratically on the distance from
$\mu$.

We conclude that the difference $i ( \lambda ' , F ( \mu ) ) - l ( \lambda ' ) $
is less than some constant times the mean of the square of the distance
from a tangent vector to $\lambda '$ to the nearest tangent vector
of $\mu$.  Let $D ( \lambda ' , \mu ) $ denote the mean square distance, and
let $d ( \lambda ' , \lambda )$ denote the distance between the
measured laminations as measured in the natural coordinates of 
$\tau$.

Let $\lambda '$ be any lamination near $\lambda$ in $\ML ( S )$
such that $\cut_{\lambda} ( \lambda ' )  =  \cut_{\lambda} ( \mu )$, and let 
$\tau$ be the train track in the chosen cover which carries $\mu$.
Any leaf of $\lambda '$ maps as a doubly infinite path on $\tau$.

The method of analysis will be to try to find for each leaf of $\lambda '$
and for each point $x$ on the leaf a leaf of $\mu$ whose trajectory on $\tau$
matches that of the leaf of $\lambda '$ as long a symmetric interval about 
$x$ as possible.  As long as the trajectories on $\tau$ match, the two leaves
remain within some uniformly bounded distance $B$ of each other in the
universal covering of $S$.  If the radius of the interval about $x$ is 
$R$, then the distance from the tangent vector to the leaf at $x$ to the leaf of
$\mu$ is not greater than $B \exp (-R )$.

How are the maximum lengths of matching trajectories distributed?
We must look a bit at what can block the existence of a matching
trajectory.  For this, it is convenient to pass to the highway picture $H_{\nu}$
associated to a train track equipped with a set of invariant weights $w_{\nu}$.
The highway is constructed as a union of Euclidean rectangles, each rectangle
corresponding to one of the branches of $\tau$.  The width
of a rectangle is the total transverse measure $w_{\nu} ( b )$,
assigned to the branch, and the length of the rectangle is 
$l ( b )$.  Each rectangle is equipped with a foliation running along its
length: the foliation is like a family of lane dividers.
The switch conditions guarantee that leaves continue from rectangle to
rectangle, except for certain exceptional leaves which branch when they arrive
at the end of a rectangle, following along the sides of the two attached
rectangles.  The lane foliation of the highway represents
the same measure class as the original lamination.

Unless $\rho  =  \tau$, the highway $H_{\lambda}$ will have some dirt roads: all
rectangles for edges of $\tau - \rho$ will have zero width, or in other words,
they are really line segments.  We may assume that
traffic entering $\rho$ on a dirt road cannot exit again on a dirt road
(without changing lanes).  This assumption is justified because there is a
fixed $\epsilon$ such that any leaf of any chain recurrent lamination
$\mu$ containing $\lambda$ which enters an $\epsilon$-neighborhood of
$\lambda$ limits on $\lambda$.  A good choice of the train track cover
will preserve this property.
\begin{figure}[htbp] \centering \approxspace{ 4in} 
\includegraphics{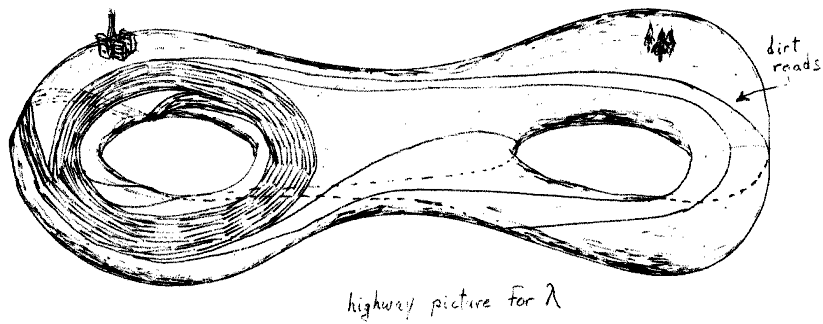}
\label{ highway picture  for lambda} 
\caption{ 
The highway corresponding to the starting lamination $\lambda$
has dirt roads (of zero width) which will be widened for certain $\lambda '$
in arbitrarily small neighborhoods of $\lambda$.  }
\end{figure}

With this proviso concerning dirt roads, there is no segment
of a leaf in the highway $H_{\lambda}$ which ends in each direction
at a divergent branchpoint in the thick part of the highway.
In other words, if you start from the middle of the segment
there has to be at least one direction to go in which the leaf
never runs into a branch point where it
splits in two, although it may run into branch points where it
unites with other leaves.  If there were any such branching, it would imply that
$\rho$ had not been chosen to have the stable topological type.

There is still another useful picture associated with $\nu$, an $\reals $
bundle $B_{\nu}$ over the train track $\tau$.
$B_{\nu}$ can be defined in terms of $H_{\nu}$ by embedding
each rectangle used in the construction of $H$ so that its two
ends lie on the two bounding lines of a strip in $\euclidean^{2}$.
Extend the gluing maps for the ends of the rectangles to gluing
maps of the bounding lines of the strips.  Topologically, the resulting space 
$B_{\nu}$ is a twisted product of $\tau $ with $ reals$.
$H_{\nu}$ is embedded in $B_{\nu}$, and the lane foliation of $H_{\nu}$
extends to a codimension one foliation of $B_{\nu}$.

$B_{\nu}$ is an example of a foliated $\euclidean^{1}$-bundle over $\tau$.
Such bundles are determined by their holonomy, which is a homorphism
from $\pi_{1} ( \tau )$ to the group of isometries of $\euclidean^{1}$,
otherwise known as an element of $H_{1} ( \tau ; O )$, where $O$
is the coefficient system of local orientations on $\tau$.
Holonomy is defined more generally for paths on $\tau$,
not just for closed loops.  The holonomy along a path is an isometry from the
fiber at one endpoint to the fiber at the other.

\begin{figure}[htbp]
\includegraphics{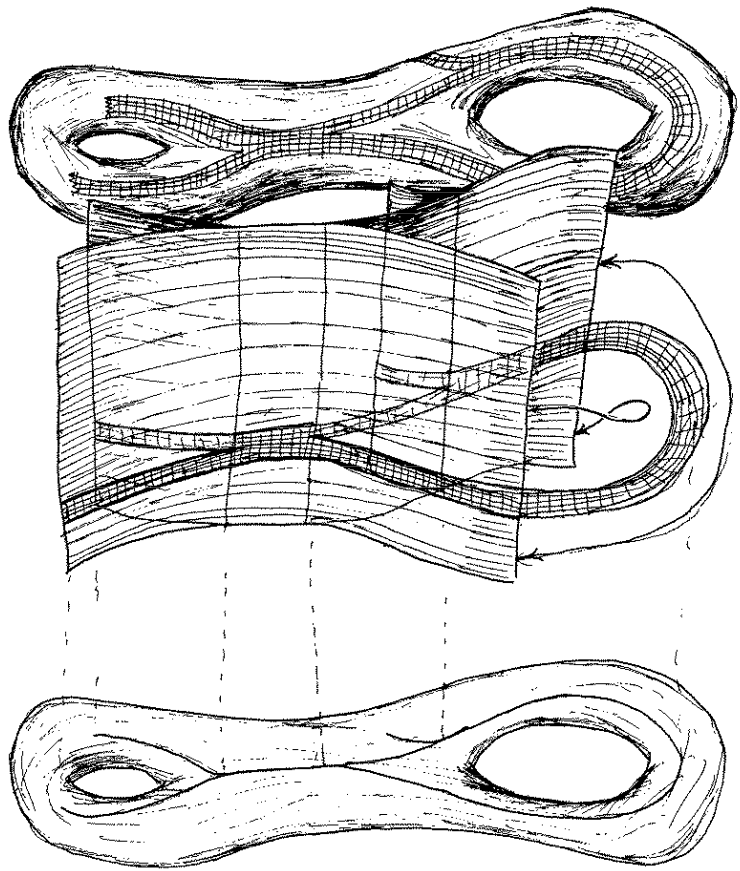}
\caption{ 
This picture shows a fragment of a train track (bottom), with
the corresponding bundle above it (middle) and the highway on the
surface behind (top).
The bundle is a branched surface, and the fibers do not have any consistent
orientation: up and down reverse above the loop on the righthand 
side of the picture.
The ramps of the highway fit neatly onto the sheets of the bundle.
}
\label{ train track and highway connected by bundle}
\end{figure}

The bundle is useful for dealing with the phenomenon that as the
highway is altered, certain journeys become impossible (without changing lanes).
The bundle tells you where you would have arrived had you been
able to take the correct exits, since it allows you to take
exits at will.

We now need to compare the two highways $H_{\lambda} '$ and $H_{\lambda}$.
Consider any point $x \element H_{\lambda} '$ not on one of the
countably many exceptional leaves.  Let $L_{R} ( x )$ be the segment of radius 
$R$ about $x$ on its leaf in the lane foliation.
If there are any trajectories of the lane foliation of $H_{\lambda}$
which match that of $x$
(that is, they follow the same path as $x$ on the track $\tau$),
the union of all such trajectories forms a
rectangle $P_{R} ( x )$ of length $2R$ and height $ h_{R} ( x )$
immersed in $H_{\lambda}$.  There is a similar rectangle 
$P '_{R} (x)$ of height $h '_{R} (  x )$ immersed in $H_{\lambda'}$
consisting of all trajectories that match that of $x$.  As $R$
increases, it probably happens from time to time that an end of
one or the other rectangle meets a divergent branch, and then
the height of the corresponding rectangle decreases.  The rectangle 
$P_{R} ( x )$ may decrease to 0 in height,
if it follows a dirt road.  Eventually, an end of the rectangle 
$P_{R} ( x )$ most likely arrives at a junction
where it is forced to take a wrong turn; the value of $R$ at that
instant is the maximum possible radius $T ( x )$ of a trajectory for $\mu$
which matches the leaf through $x$.

Draw a line segment $\beta$ transverse to the highway
which extend the end of the rectangle $P_{T} ( x )$
to the branch point $p_{1} \in H_\lambda$
where the two roads diverge.  This end of the rectangle
is in the thick part of $H_{\lambda}$: since
$\cut_{\lambda} ( \lambda ') = \cut_{\lambda} ( \mu )$, no part of $P_{R}$
on a dirt road of $H_{\lambda}$ is ever forced to take a wrong turn.

Tracing the leaves of the lane foliation of $H_{\lambda}$ backwards along
the side of $P_T$ from $\beta$,
one finds that all the lanes through $\beta$ have to peel away from
$P_T$ before they reach the other end of $P_{T} $.  Let $p_{2}$
be the last point where lanes peel away;
it is a branch point of $H_\lambda$ somewhere along the
edge of $P_{T}$.
It may happen that $p_{2}$ is on a junction of a wide road with a dirt road,
but it is not on a dirt road.

\begin{figure}[htbp] \centering \approxspace{5in} 
\includegraphics{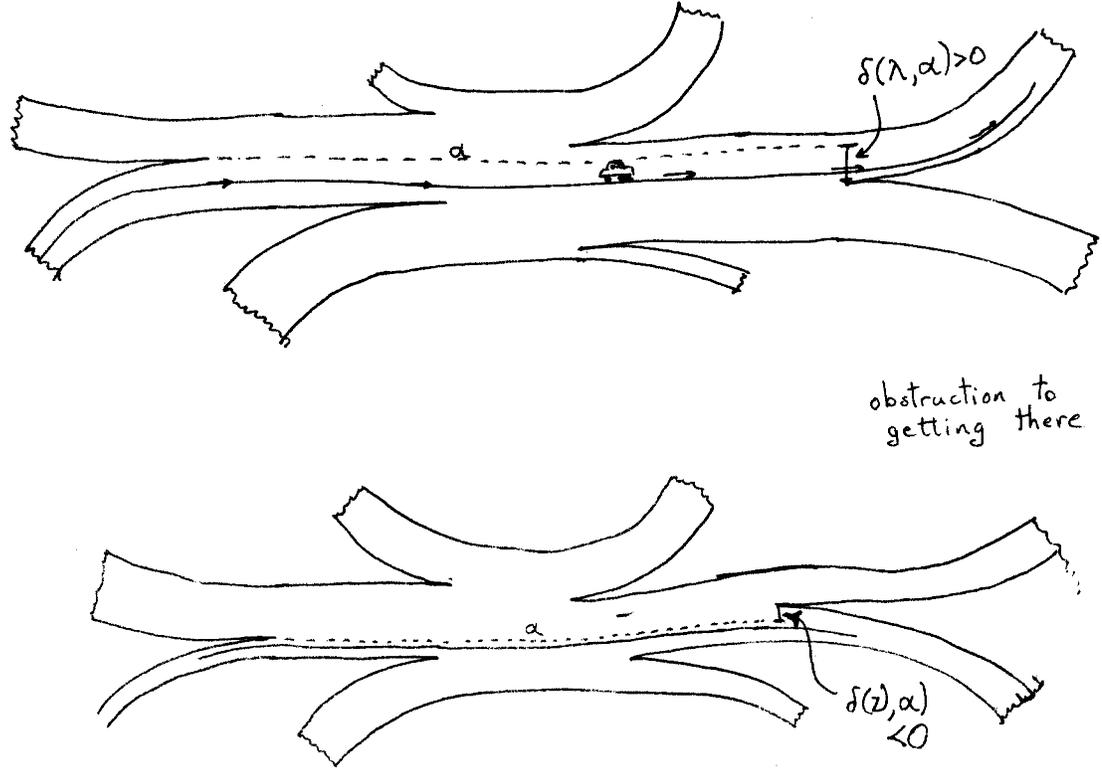}
\label{ obstruction to getting there} 
\caption{
A homotopy class of arcs between two divergent branch points
represents the obstruction to existence of a longer matching trajectory.
}
\end{figure}
Let $ \alpha = \alpha ( x ) $ denote the homotopy class of paths joining 
$p_{2}$ with $p_{1}$, following the side of $P_{T}$ and then $\beta$.
We think of $\alpha ( x )$ as representing the obstruction to finding
any longer trajectory on $H_{\lambda}$ matching the trajectory of 
$x$.  For any measured lamination $\nu$ carried by $\tau$, there is associated
to the homotopy class $\alpha$ a number
$$ \delta ( \nu , \alpha )  = \hol_\alpha ( p_2 ) - p_1
$$
which measures the signed distance between the lane of $p_{2}$ and
that of $p_{1}$
at the time it reaches that junction.  The subtraction here makes sense, since 
$\hol_{\alpha}$ maps the fiber of the bundle over $p_{2}$ to the fiber over 
$p_{1}$.  The sign is determined once we choose an orientation for the fiber
at $p_{1}$.  From the bundle picture $B_{\nu}$ it follows that
$\delta ( \nu , \alpha )$ is defined whether or not the lane of 
$p_{2}$ ever in fact reaches the junction.
It is a linear function in train track coordinates. 

In the present case, we see that $\delta ( \lambda , \alpha )$ and
$\delta ( \lambda ' , \alpha )$ have opposite signs, and that
$| \delta ( \lambda ' , \alpha ) | \ge h '_{T}$, while
$| \delta ( \lambda , \alpha ) |$ is the length of the arc $\beta$.
Note that the length of the path $\alpha$ on $\tau$ is not greater
than $2 T$.

As the point $x$ varies over the highway $\lambda '$, the obstructing path 
$\alpha ( x )$ will vary.  Any path which might occur as $\alpha ( x )$
is represented as a simple path on the highway $H_{\lambda}$. 
Just as a simple closed curve or a measured lamination, it determines
a set of weights on the branches of $\tau$.  The weights are all integers,
and they satisfy the switch conditions
except at the switches at the two ends of the path, where they
fail the switch conditions by a fixed amount: at these
two switches, if the two converging branches
have weights $a$ and $b$ and the third branch has weight $c$, then
$ a + b + 1 = c.$
From this information, the local picture near the switch is easily
reconstructed.

We conclude that the number $N ( R )$ of possible paths $\alpha$
is dominated by a polynomial of $R$ with degree
at most the dimension of the measured lamination space (depending on
how big $\tau$ is).

Let us now estimate the mean square distance $D ( \lambda ' , \lambda )$
from a tangent vector of $\lambda '$ to the closest tangent vector of
$\lambda$.  We take the mean with respect to the product of Lebesgue measure
along the leaves with the transverse invariant measure.
This measure is about the same as two-dimensional Lebesgue measure on
$H_{\lambda} '$.  For any point $x \element H_{\lambda} '$,
the distance from the tangent vector at its image in $\lambda '$
to the nearest tangent vector in $\lambda$ is less than some constant
times $\exp ( - T ( x ) )$.  The width of the rectangle of trajectories
associated with $x$ is not greater than $\delta ( \lambda ' , \alpha ( x )$.
The length $2 T(x)$ of the rectangle $P_T(x) (x)$ is
at least as great as the length of $\alpha (x)$,
but there probably are also longer and
shorter rectangles $P_T(y) (y)$ having the same obstructing arc
$\alpha(y) = \alpha(x)$.
Integrating over all $y$ with $\alpha ( y ) = \alpha ( x )$,
we find that the net contribution to $D ( \lambda ' ,  \lambda )$
is less than a constant (depending on $\tau $) times
$\delta(\lambda' ,\alpha) \exp(-\length (\alpha))$.

Keeping in mind that $\delta(\lambda', \alpha)$ and $\delta(\lambda, \alpha)$
are of opposite sign, we see that
\begin{align*}
D ( \lambda ' , \lambda ) & < C  \sum_{\alpha }
\delta ( \lambda ' , \alpha ) \exp ( - {\length ( \alpha )} ) .
\\
&< C  \sum_{\alpha }
| \delta ( \lambda ' , \alpha )  - \delta ( \lambda , \alpha ) | \,
\exp ( - {  \length ( \alpha ) } ) .
\end{align*}

The holonomy  $\delta ( * , \alpha )$ is a linear function in train track
coordinates.  If expressed in terms of a basis for the linear functions,
the total absolute value of the coefficients is less than some constant times
the number of branches which $\alpha$ crosses, which is less than
some constant times the length of $\alpha$.  Therefore, the difference
$| \delta ( \lambda ' , \alpha )  - \delta ( \lambda , \alpha ) |$
is dominated by some constant times the length of $\alpha$ times the distance 
$d ( \lambda ' , \lambda )$ in measured lamination space.
Grouping the contributions from different arcs
according to the intervals which the lengths lie in, we get an estimate
$$ D ( \lambda', \lambda )
< K  \sum_1^\infty N' ( R )  R  \exp ( -R ) d ( \lambda', \lambda ) .
$$
We can improve this estimate by noting that when  $\lambda '$
is close enough to $\lambda$, any particular arc $\alpha$ stops
acting as an obstruction.  This follows from the fact that 
$\delta ( \lambda ,  \alpha )$ is never zero, so that 
$\delta ( \lambda ' , \alpha )$ will eventually have the same sign.
Therefore, for all $N > 0$, for $\lambda '$ sufficiently close to 
$\lambda$ we will have
$$ \frac {D ( \lambda ' , \lambda )}{ d ( \lambda ' , \lambda ) }
< K  \sum_N^\infty Q ' ( R )  R  \exp ( -R ) .
$$

Since the left hand side tends to 0 as $N$ tends to infinity, this implies
that the affine function $i ( \lambda ' , F ( \mu ) )$
gives the derivative of the length function.

Next we must check continuity of the formula over the entire tangent cone 
$C ( \lambda )$, or equivalently, the continuity of $i ( \lambda ' , F ( \mu )$
in a neighborhood of $\lambda$.  The formula depends on the geodesic lamination 
$ \mu $, which we should write $ \mu_{\lambda '}$ to emphasize dependence on 
$\lambda '$, or at least $\cut_{\lambda'}$.  For a measured lamination
$\lambda'' $ near $\lambda ' $, most of its measure is near the support of 
$\lambda '$.  The foliations $F ( \mu_{\lambda '} )$ and
$F ( \mu_{\lambda''} )$ nearly agree in this neighborhood, even though 
$ \mu_{\lambda''} $ may have new strands outside of this neighborhood.
This easily implies continuity of the first-order approximation to length,
$$ i ( \lambda', F ( \mu_{\lambda'}) )  ,
$$
and hence continuity of the derivative (as a function of tangent vectors
at one point only).

The error estimate for the first-order approximation to length depended on the
numbers $\delta ( \lambda , \alpha )$ and on the length of the arcs $\alpha $
in the hyperbolic metric, and constants concerning the geometry of the track
in the surface.  Clearly, when $\lambda$ is fixed, the estimate depends
continuously on the geometry of the surface.
\endproof %{\ref {Derivative is intersection number}}

\section {Maximal ratio laminations} \label {Maximal ratio laminations}

The next theorem will (among other things) imply a kind of local
uniqueness for the lamination which attains the supremum of $r$.

\begin{theorem}[Parametrized conjugacy of laminations]\label{Parametrized conjugacy of laminations}
Let $\lambda$
be any measured lamination (up to homotopy) on a closed surface $S$,
let $g$ and $h$ be two hyperbolic structures on $S$, and let $\lambda_g$
and $\lambda_h$ denote the realizations of $\lambda$ as geodesic measured
laminations in the two metrics..
Then if the ratio $\length_{h}(\lambda_h)/\length_{g}(\lambda_g)$ is maximal
among all laminations
whose support is contained in some neighborhood of the support of $\lambda$,
there is a homeomorphism from a neighborhood of $\lambda_{g}$ to a neighborhood
of $\lambda_{h}$, in the isotopy class of the identity,
which takes each leaf of $\lambda_{g}$ linearly to its image leaf and whose
Lipschitz constant equal to the ratio of lengths.
\end{theorem}

Note that every simple closed geodesic trivially satisfies the hypothesis
of this theorem. It also trivially satisfies the conclusion.

The realizations of laminations on different surfaces
are always related by a homeomorphism of the surface: what is special
is that the homeomorphism is only a linear change of parametrization.

\proof %{Proof of \shortlabel Parametrized conjugacy of laminations:8.1:}
Under the hypotheses of the theorem, it follows that the derivative
of the ratio of lengths $ r_{gh}$ must be zero at least in
the near subspace of the tangent space.

This means that the measured foliations $F_{g} ( \lambda )$
must have the measure class which is the multiple $r_{gh} ( \lambda )$
of the measure class of $F_{h} ( \lambda )$.

Map leaves of $F_{g} ( \lambda )$ to corresponding leaves of $F_{h} ( \lambda )$
while taking leaves of $\lambda$ to corresponding leaves of $\lambda$.
This gives the desired Lipschitz map, when restricted to a sufficiently
small neighborhood of $\lambda$
(see \ref{Neighborhood of a lamination}.)
\endproof %{\ref {Parametrized conjugacy of laminations}}

This theorem suggests that the notion of laminations which
maximize the ratio of length between two surfaces can be extended to more
general geodesic laminations, not just measured laminations.
We can say that a geodesic lamination $\lambda$ is {\it ratio-maximizing}
for a pair of hyperbolic structures $g$ and $h$ if there is a $r$-Lipschitz map
from a neighborhood of $\lambda$ in $g$ to a neighborhood of $\lambda$ in $h$,
in the correct homotopy class, where $r$ is the supremum of $r_{gh} ( \mu )$
for measured laminations.

\begin{theorem}[mu(g,h) well-defined]\label{mu(g,h) well-defined}
For any two distinct hyperbolic structures $g$ and $h$,
there is a unique maximal ratio-maximizing chain recurrent lamination
$\mu ( g , h )$ which contains all other ratio-maximizing chain recurrent
laminations.
\end{theorem}

\proof %{Proof of \shortlabel mu(g,h) well-defined:8.2:}
Consider any two chain-recurrent laminations $\lambda_{1}$ and $\lambda_{2}$
which are ratio-maximizing.  Suppose that they cross somewhere at an angle;
we will show that from them we can construct a curve with a larger value of 
$r_{gh}$, giving rise to a contradiction.

The first step is show

\begin{lemma}[approximate chain recurrent by curve]\label{approximate chain recurrent by curve}
For any chain recurrent lamination $\lambda$ there is a closed
curve $\gamma$ which approximates $\lambda$ arbitrarily closely in the
Hausdorff topology and whose length is within an arbitrarily
small additive constant of the intersection number with $F ( \lambda )$.
\end{lemma}

\proof[Proof of \ref{approximate chain recurrent by curve}]
It is convenient to pass to the oriented double cover $\mu$ of $\lambda$, unless
$\lambda$ is already orientable, in which case we pick an orientation and
call it $\mu$.  The $\alpha$ limit set of a leaf in $\mu$ is defined
as limit set of the leaf in negative time; the $\omega$ limit set is
the limit set in positive time.

There is a directed graph whose vertices are leaves of $\mu$ and whose
directed edges join pairs of leaves such that the $\omega$ limit set of the first
intersects the $\alpha$ limit set of the second.  The $\alpha$ and $\omega$
limit sets of a leaf are closed invariant subsets of the lamination.
There are only a finite number of possible candidates for the $\alpha$
limit set of a leaf, and only a finite number of possible candidates for the 
$\omega$ limit set.
There are also only a finite number of isolated leaves which are not in any
$\alpha$ or $\omega$ of any leaf limit set.
Consequently, if we form an equivalence relation on the vertices
of the graph that two leaves are equivalent if each is in both the $\alpha$
and $\omega$ limit set of the other (so that in particular, they have the same 
$\alpha$ and $\omega$ limit sets), we obtain a finite directed graph.

The hypothesis that $\lambda$ is chain recurrent implies that this finite
directed graph admits a closed loop through each edge.  The hypothesis that 
$\lambda$ is connected implies that there is a single closed loop which
passes through every edge.  Given such a closed loop, a closed curve on 
$S$ can be constructed by first picking a leaf of $\mu$
representing each vertex in the loop.  The curve $\gamma$ on $S$ follows a
long segment in one of the leaves until it comes within $\epsilon$
of the next leaf, then hops over, and repeats the process until it closes.
\endproof %[\ref {approximate chain recurrent by curve}]

\proof [Continuation of proof of \ref {mu(g,h) well-defined}]
Apply the lemma to the $\lambda_{i} $, finding curves $\gamma_{i}$
which approximate the $\lambda_{i}$ well for both metrics.

When two long curves meet, the angle between them is ambiguous:
either of two complementary angles $\theta$ and $\pi - \theta$ is correct.
The two angles correspond to two different ways of cutting and
regluing to yield a single curve.  There is a function 
$f ( \theta)$, which decreases monotonically from $\infty$ to $0$ as $\theta$
goes from $0$ to $\pi$, such that when two very long geodesics are cut and
reglued at an angle $\theta$, the geodesic in the new homotopy class is
the sum of the two lengths minus $f ( \theta )$, up to a small correction
which for fixed $\theta$ decreases exponentially as the geodesics lengthen.

Pick one of the angles formed by the $\gamma_{i}$ which does not
decrease going from metric $g$ to metric $h$.  Cut and reglue.
After applying the various estimates of length to the resulting curve $\gamma$,
it is clear that the value of $r_{gh} ( \gamma )$ is larger than that for 
$\lambda_{1}$ and $\lambda_{2}$.
\endproof %{\ref {mu(g,h) well-defined}}

The next theorem shows that $\mu ( g , h )$ behaves nicely as $g$ and $h$
vary over Teichm\"uller space.

\begin{theorem}[local behavior of mu(g,h)]\label{local behavior of mu(g,h)}
Let $g$ and $h$ be any two distinct hyperbolic structures on $S$.  If $g_{i}$
and $h_{i}$ are sequences of hyperbolic structures converging to $g$ and $h$
respectively, then $\mu ( g , h )$ contains any lamination in the limit set of 
$\mu ( g_{i} , h_{i} )$ in the Hausdorff topology.
\end{theorem}

\proof %{Proof of \shortlabel local behavior of mu(g,h):8.4:}
The idea is to make use of the Lipschitz maps defined on
neighborhoods of the $\mu ( g_{i} , h_{i} )$.
It would be nice if we could pass to a limit.
The difficulty is that these neighborhoods do not have a uniform size.

The problem begins with the definition of $F ( \mu )$, which depends on picking
a neighborhood of $\mu$ which has the stable topological type.
There is a related definition for a transverse foliation $F_{\epsilon} ( \mu )$
which is defined in an entire $\epsilon $ neighborhood of $\mu$, where
$\epsilon$ can be chosen once and for all, independently of $\mu$.
First $F_{\epsilon} ( \mu )$ is defined in the $\epsilon$-thin part of
the cut surface $S_{\mu}$, with leaves being either horocycles or 
equidistant curves from the core geodesic, as appropriate.
Then it is extended to the epsilon neighborhood of the boundary of $S_{\mu}$
by compromising between the different definitions, but always perpendicular
to the boundary.  It is extended to the entire $\epsilon$-neighborhood of 
$\mu$ by integrating.
The new definition agrees with the old definition, except that it is defined
in a larger neighborhood.

The lamination $\mu ( g , h )$ has the propety that the measure class of 
$F_{h} ( \mu ( g , h ) )$ is a multiple of the measure class of
$F_{g} ( \mu ( g , h ) )$.
It need not be true that the measure class of
$F_{{h,} \epsilon} ( \mu ( g , h ) )$ is a multiple of that of
$F_{{g,} \epsilon} ( \mu ( g , h ) )$: the equidistant curves which cut
across a thin part of $S_{\mu}$ might change by a small transverse shear.
In the first place, the topology of the $\epsilon$-thin part may change
between the two metrics.
The easiest way to get around this problem is to note that the measure
class of the foliation depends only on the geometry of the surface, the
lamination, plus the topology of the neighborhood in which it is defined.
We think of $F_{{\epsilon} ,g} ( \mu )$ as being defined in some
neighborhood of $\mu$ union a finite set $\Delta$ of disjoint
geodesic arcs perpendicular to two boundary components of the cut surface 
$S_{\mu}$.  For any such set $\Delta$, there is a well-defined measure class
of foliations $F_{{\Delta} , g} ( \mu )$.  Now it is possible to compare 
$F_{{\Delta} , g } ( \mu (g , h)$ with $F_{{\Delta} , h } ( \mu ( g, h)$.

Let $\delta \element \Delta$ be any short geodesic perpendicular to two
boundary components of $S_{{\mu} ( g , h )}$ in the homotopy class.
Associated with $\delta$ is an error $e ( \delta )$, defined as the
amount of shear transverse to $\delta$ required to change
$\exp ( K ( g , h ) ) F_{{\Delta} , g } ( \mu ( g , h ) )$ to
$F_{{\Delta} ' , h } (\mu ( g , h ) )$.  Actually, $e ( \delta )$
is well-defined only if $\mu ( g , h )$ is non-orientable.  The reason is that
in the non-orientable case, the map taking
$\exp ( K ( g , h ) ) F_{g} ( \mu ( g , h )$ to $F_{h} ( \mu ( g , h )$
is unique, up to restriction to a smaller
domain and an isotopy which takes each leaf to itself.
In the orientable case, there is a whole one-parameter family of different
maps, generated by a flow coming from a vector-field extending the tangents
to $\mu ( g , h )$.

We will come back to the orientable case later.  In the non-orientable case, we
claim that there is a upper bound to $e ( \delta )$ which depends only on
$K ( g , h )$ and on the length of $\delta$ for $g$.  To show this,
use the proof of lemma \ref{approximate chain recurrent by curve}
to find a map of an interval $\gamma$, closely approximating 
$\mu ( g , h )$ in the sense of the lemma, and connecting one end of 
$\delta $ to the other end in the direction such that the shear will tend to
increase the length of $\delta \union \gamma$ in the metric $h$.
If the geodesic arc $\delta$ happens to be longer in the metric $h$,
this will also tend to increase the length of $\delta \union \gamma$.
If $\delta$ is shorter in $h$ , this will tend to decrease $\delta \cup\gamma$
in $h$, but only by an amount which goes rapidly to zero with the length
of $\delta$ in $g$.  The fact that the new curve cannot have a value of 
$r_{gh}$ greater than $r_{gh} (\mu ( g , h ))$ gives the desired bound
for $e ( \delta )$.

In the orientable case, if there is any $\delta \element \Delta$
which joins two like-oriented leaves, then $e ( \delta )$ is actually
well-defined, and the same argument applies.  If $\Delta$ is empty,
there is nothing to prove.

Otherwise, all elements of $\Delta$ join leaves of which point in opposite
directions, and $e$ is defined up to addition of a constant If $e$
is a constant, then the map taking $\exp ( K ( g , h ) ) F_{g} ( \mu ( g , h )$
to $F_{h} ( \mu ( g , h )$ can be adjusted to extend to $F_{\Delta} $.
If $e$ is not a constant, then a curve on $S$ can be formed using two elements
of $\Delta$ for which $e$ takes different values, together with two intervals
mapped into $S$ joining their endpoints which closely approximate
$\mu ( g , h )$.

Now we can return to the main proof.  Suppose that $\{ g_{i} \} \arrow g$
and $\{ h_i \}\arrow  h$, and that $\mu ( g_{i} , h_{i} )$ has a Hausdorff
limit $\mu$.  Let $U_{g}$ and $U_{h}$ be small neighborhoods of $\mu$
of the stable topological type.  We can assume that 
$\mu_{g} ( g_{i} , h_{i} )$ is contained in $U_{g}$ and comes within 
$\epsilon$ of every point in $U_{g}$, and we can let $\Delta_{i}$
be the set of arcs necessary to augment a regular neighborhood of 
$\mu_{g} ( g_{i} , h_{i} )$ to be the same homotopy class as $U$.
Elements of $\Delta_{i}$ are arbitrarily short, as $i\arrow \infty$.
From the discussion above, it follows that the measure class of
$F_{h} ( \mu ) = \exp ( K ( g , h ) ) F_{g} ( \mu )$.
The theorem follows from an application of
\ref{mu(g,h) well-defined}.
\endproof %{\ref {local behavior of mu(g,h)}}

This theorem can be interpreted as saying that $\mu (g , h )$ is
continuous in a non-Hausdorff topology on the set of chain recurrent
laminations, where a neighborhood of a lamination consists of all
laminations contained in a neighborhood of the lamination on the surface.

Here is the theorem which originally motivated this paper:
\begin{theorem}[constructing stretch geodesics]\label{constructing stretch geodesics}
For any two hyperbolic structures $g$ and $h$,
the minimum Lipschitz constant of a homeomorphism
is the same as the maximum ratio of
lengths of laminations, that is,
$$ K(g, h ) = L ( g , h ) .$$
There is a geodesic with respect to these metrics consisting of
a finite concatenation of stretch paths, each of which stretches
along some lamination containing $\mu ( g , h )$.
\end{theorem}
\proof %{Proof of \shortlabel constructing stretch geodesics:8.5:}
Given $g$ and $h$, let $\mu ( g , h )$ be the maximal chain recurrent
lamination maximizing stretch between $g$ and $h$.
Choose arbitrarily a maximal lamination $\mu$ which extends 
$\mu ( g , h )$.  Consider the stretched surfaces $g_{t}$ using $\mu$.
By theorem \ref{local behavior of mu(g,h)}, the lamination 
$\mu ( g_{t} , h )$ must be contained in a small neighborhood of 
$\mu ( g , h )$ for $t$ near 0.
By construction of the stretch maps, there is a Lipschitz map
defined on a small neighborhood of $\mu$ with Lipschitz constant 
$\exp( K ( g , h ) -t )$ from the surface at time $t$ to $h$.
Therefore, the only maximally-stretched
laminations contained in this neighborhood are contained in 
$\mu $.

If $\mu ( g_{t} , h )$ is not constant until $K ( g_{t} , h ) = 0$,
then there is some first $t_{1}$ for which
$\mu ( g_{t_1} , h ) \ne  \mu ( g , h )$.  Again by
\ref{local behavior of mu(g,h)}, $\mu ( g_{t_1} , h ) $ contains 
$\mu ( g , h )$.

Repeat the process.
The maximal maximally-stretched lamination can only increase a finite number
of times, since there is a bound to the length of an increasing sequence
of laminations on any given surface.
Therefore, the process continues until it reaches a surface 
$g '$ such that $K ( g ' , h ) = 0$.  By theorem \ref{K(g,h) positive},
$g ' = h$.
\endproof %{\ref {constructing stretch geodesics}}

\section {Cataclysm coordinates for Teichmuller space} \label {Cataclysm coordinates for Teichmuller space}

The construction of proposition \ref{Neighborhood of a lamination}
points to a family of global
coordinate systems for Teichm\"uller space which are analogous to
the upper half-space picture for hyperbolic space:
A hyperbolic surface can be described in terms of the geometry
in a neighborhood of some maximal lamination $\mu$.

To make this into an effective set of coordinates, we have two tasks.
First, we will give a better description of the space of
transverse measured foliations defined in a neighborhood of $\mu$.
Second, we will show that when $\mu$ is maximal, there is a unique choice of
sharpness functions for any transverse measured foliation
so that the hyperbolic structure defined in a neighborhood of $\mu$
extends to a complete hyperbolic structure on all of $S$.

Let $\mu$ be a maximal geodesic lamination on $S$.
Consider any measured foliation $F$ defined in a neighborhood of $\mu$
which is transverse to $\mu$.  For each ideal triangle of $S_{\mu}$, 
$F$ restricts to a foliation defined in all but a compact subset of the
interior.  The claim is that there is an extension of $F$ to the entire triangle
with one three-pronged singularity in the interior, which is unique up
to an isotopy supported on a compact subset of the interior.
In fact, constructing the foliation is a matter of solving a system of three 
linear equations, where the unknowns represent the total additional
transverse measure of leaves passing between a pair of sides.
In conclusion, we have

\begin{proposition}[Global F(mu)]\label{Global F(mu)}
For any hyperbolic surface and any maximal geodesic lamination, the
measured foliation $F ( \mu )$ extends uniquely to a measured
foliation (also denoted $F ( \mu )$) on all of $S$.
\end{proposition}

A related fact is that for a maximal geodesic lamination $\mu$,
no other data than $F ( \mu )$ is needed to determined a hyperbolic
structure on all of $S$.

\begin{proposition}[unique extendible neighborhood]\label{unique extendible neighborhood}
If $\mu$ is any maximal lamination on a complete hyperbolic surface
$S$ and if $F$ is any measured foliation transverse to $\mu$,
defined in a neighborhood of $\mu$, and standard near any cusps of $S$,
then there is a unique complete hyperbolic structure on $S$ such that
$$ F = F ( mu )
$$
\end{proposition}
\proof %{Proof of \shortlabel unique extendible neighborhood:9.2:}

The components of the cut surface $S_{\mu}$ are all ideal triangles.
If the hyperbolic structure obtained from proposition
\ref{Neighborhood of a lamination}
extends to each of these triangles, then it gives a complete hyperbolic
structure on $S$; conversely, it is certainly necessary that the
structure defined on the boundary of each ideal triangle extend to a triangle.

In an ideal triangle, there is a canonical measure class of foliations
each of whose leaves is isotopic rel endpoints to a horocycle.
(The horocycles themselves leave a little gap in the middle, so they
don't make a measured foliation as is).  Isomorphisms from $F$ restricted to
the components of $S_{\mu} $ to the canonical foliation of an ideal triangle
determine the sharpness functions.
\endproof %{\ref {unique extendible neighborhood}}

Any train track approximation of $\mu$ determines a coordinate system for a
subset of foliations transverse to $\mu$, as described in section
\ref{Derivatives of lengths of laminations}
A transverse foliation is described by a positive function on the branches
of the track, up to addition of linear combinations of the switch vectors.
In the case of cusps, we can assume that the train track approximation
consists of one branch going out to infinity, and we
assign infinity to this branch.
These functions, up to the equivalence relation,
form a convex cone in a certain vector space.

As the train track approximation gets finer and finer, the train
tracks can be made transverse to more and more
of the measured foliations transverse to 
$\mu$, so that the coordinate domains get larger and larger.
There are natural injective linear
maps from the coarser transverse coordinate neighborhoods
to the finer ones.  Thus, the set of all measured foliations transverse to $\mu$
and standard near cusps can be described as an increasing union of
convex cones, all of the same dimension.  The union is itself a convex cone.
What does it look like?

\begin{theorem}[linear characterization of transverse foliations]\label{linear characterization of transverse foliations}
Let $\tau$ be any train track carrying $\mu$.  The set of measured foliations
transverse to $\mu$ is linearly isomorphic to the set of functions on the
branches of $\tau$ assigning infinity to each noncompact branch of 
$\tau$ which when paired with any transverse invariant measure for $\mu$
are positive, modulo the relations.
\end{theorem}

The set of transverse invariant measures for $\mu$ forms the cone
on a simplex of some dimension (as for any foliation):
every transverse invariant
measure decomposes in a unique way as a convex combination of ergodic measures.
The convex cone above is the dual of the cone on some simplex.
The simplest case is when $\mu$ has a unique transverse invariant measure;
then the set of measured foliations transverse to $\mu$ is parametrized
by a half space.

This gives a coordinate system for Teichm\"uller space which
is analogous to a certain coordinate system for hyperbolic space,
obtained by a parallel projection of the Lorenz model along a family of 
rays parallel to one of the light rays.
The image is a half-space, but the model is not like the Poincar\'e half-space
model.  All the finite boundary points in this model are obtained from
one finite boundary point in the Poincar\'e disk model, while the parallelism
classes of rays going to infinity all converge to different points in
the Poincar\'e disk model.

\proof %{Proof of \shortlabel linear characterization of transverse foliations:9.3:}
For segment $I$ of a leaf of $\mu$, the image on $\tau$ passes
through a sequence of branches of $\tau$.  Let $S(I)$
be the sum of the numbers associated to this sequence of
branches, weighted proportionally to the length for the two end
branches may be only partially covered by $I$.
Length is measured with respect to some metric on $\tau$.  If $I$
intersects a noncompact branch of $\tau$ (which is assigned infinity),
we can make the convention that the weight is the length of the interval
of $I$ on this branch.

The key point is that under the hypotheses, there is a number $L$ such
that for any interval of length $I$ of length greater than $L$,
$S ( I ) > 0$. It suffices to prove this fact for the case that $I$
does not intersect any noncompact branch.  The proof is by contradiction.
If there were no such number, then there would be a sequence $I_{k} $
of segments whose lengths go to infinity, such that $S ( I_{k} ) < 0$,
the length of $I_{k} )$ goes to infinity, and $I_{k}$ maps
to the compact branches of $\tau$.  But the ``hitting measure'' for 
$I_{k}$ , normalized by dividing by the length of $I_{k}$, is almost invariant;
a subsequence of these measures converges to a transverse invariant measure for 
$\mu$ . Its pairing with the weighting of $\tau$ is nonpositive.

Once a number $L$ is found, we can directly construct
a transverse measured foliation as follows.
First, construct a neighborhood of $\mu$ with two transverse foliations,
the foliation $G$ which extends $\mu$, and the foliation $F$ which is
transverse to $G$ and whose leaves map to points on $\tau$.
There is a subneighborhood $U$ of $\mu$ consisting of points whose
$G$-leaf contains a segment of radius $L$, as measured by the
map to $\tau$.

We shall construct a transverse invariant measure for $F$ restricted
to $U$.  The idea is that the original set of weights for $\tau$ define
a transverse invariant signed measure for $\tau$ on the larger neighborhood.
One can think of this as a closed 1-form with coefficients in the
$\reals$-bundle associated with local transverse orientations.
Averaging over long stretches will produce a non-singular closed 1-form.
Explicitly, let $\sigma$ represent the transverse measure for $F$ induced
from arc length on $\tau$.  For any point $p$ in $U$, let 
$T ( p ) = S ( I ) / 2L$, where $I$ is the $G$-leaf of $p$ of radius $L$.

$I$ also follows the trajectory of some $\mu$-leaf, so $T ( p )$
is positive.  Then $T_* \sigma$ defines a new transverse invariant measure for 
$F$, which is what we want.

The construction is independent of the choice of $L$, up to measure
equivalence and restriction to a smaller neighborhood.
\endproof %{\ref {linear characterization of transverse foliations}}

\begin{proposition}[transverse laminations and transverse foliations]\label{transverse laminations and transverse foliations}
Let $\mu$ be a maximal geodesic lamination.  Foliations such as $F$
above are in one-to-one correspondence with compactly supported measured
laminations $\lambda$ on $S$ such that each leaf of $\lambda$ intersects 
$\mu$ transversely infinitely often, and each leaf of $\mu$ which does not
go to a cusp intersects $\lambda$ transversely infinitely often.
(In counting intersections, we think of leaves as parametrized by $reals$).

\end{proposition}
\proof %{Proof of \shortlabel transverse laminations and transverse foliations:9.4:}
There is a standard construction to derive a measured lamination from a
measured foliation.  In this case, we first throw away all leaves of 
$F$ concentric to the cusps; each remaining leaf is homotopic to a geodesic.
The associated lamination $\lambda$ consists of the closure
of the set of geodesics obtained in this way, and it has
a transverse invariant measure induced from the measure for $F$.

It is clear that each leaf of $\lambda$ intersects $\mu$ transversely infinitely
often.  A leaf of $\lambda$ intersects leaves of $F$ infinitely often;
if the leaf does not go out to a cusp, the leaves it intersects cannot
be the ones that have been thrown away, so they go over to leaves of $\lambda$.

One way to reconstruct a measured foliation from $\lambda$ is to use
a fine train track approximation of $\mu$, fine enough that $\tau$ is
transverse to $\mu$.  To each compact branch of $\tau$, assign the total measure
of intersection with $\lambda$, and to each noncompact branch assign infinity.
\endproof %{\ref {transverse laminations and transverse foliations}}

A lamination $\lambda$ satisfying the hypothesis of the proposition
is {\it totally transverse} to $\mu$.

Putting everything together, we have a parametrization of Teichm\"uller
space by a convex cone in affine space,
the dual cone to the space of transverse invariant measures
supported on $\mu$.  The parametrization depends on $\mu$.
One case of special interest is when $S$ has cusps and $\mu$
consists of a finite number of leaves going out at both ends to the cusps:
that is, $\mu$ consists of the edges of an ideal triangulation.
There are no compactly-supported measures on $\mu$, so the parameter
space is the dual of the empty set, that is, an entire vector space.
The parameter describes how the triangles are to be glued together.
For the interpretation of the parameter as a compactly-supported measured
lamination,
the case of the empty lamination corresponds to the symmetric gluing of
the triangles.

The parametrization is real analytic, for the same reason as with
the stretch maps:
the holonomy around any element of the fundamental group is expressed
as an infinite product of elements in a nonabelian group
which converges well.
The terms in the product depend in a simple analytic way on the linear
coordinates, so the product is real analytic.
(In fact, the product converges uniformly for a neighborhood of the parameter
space in its complexification).

In view of theorem
\ref{Derivative is intersection number}
the parameter can also be interpreted as
the first derivative of lengths of laminations at some measure
of maximal support on $\mu$ along a linear fragment determined by $\mu$,
when $S$ is compact and $\mu$ is chain recurrent.

\begin{corollary}[Analyticity of the derivative]\label{Analyticity of the derivative}
The first derivative of length of laminations along any linear fragment
of the tangent space to measured lamination space
is a real analytic function of the hyperbolic structure.
\end{corollary}

\proof %{Proof of \shortlabel Analyticity of the derivative:9.5:}
It is not so easy to see directly that the formula for the 
derivative is analytic,
but it is immediate from the analytic parametrization of Teichm\"uller
space (which is an analytic diffeomorphism in the opposite direction).
\endproof %{\ref {Analyticity of the derivative}}

There is another form for the parameters, the shear coordinates,
which can be obtained by a linear
transformation from the transverse train track coordinates.
The shear coordinates give another way to describe measured foliations
transverse to a maximal geodesic lamination; unlike the transverse coordinates,
they do not work to describe foliations in a neighborhood when the
lamination is not maximal.
The coordinates are defined by a close train track approximation $\tau$ of $\mu$
which is transverse to $F$.  Then for any branch $b$ of $\tau$, there is a
triangle of $S - \tau$.  Inside each of these triangles, there is one
3-pronged singularity of $F$, one of whose prongs hits the side of the
triangle containing $b$.  Assign to $b$
the real number whose absolute value is the total
transverse measure between the prongs from either side, and whose
sign is positive or negative depending on
whether the prongs are displaced by a left shear or a right shear.
(This is well-defined if $S$ has an orientation; otherwise,
a sign convention can be separately chosen for each $b$.)

The linear formula for the shear coordinates in terms of the
transverse coordinates is fun to work out.
It turns out that the total shear across a branch is
one-half of the alternating
sum of the total transverse measures assigned to the sides of
the polygon obtained by gluing the two triangles together.
This polygon may be either a quadrilateral, a hexagon (with
one inward-pointing vertex) or an octagon (with two inward-pointing
vertices).
Note that the formula is not changed by adding switch vectors to
the transverse weights.
The shear coordinates in the oriented case satisfy the switch relations,
so unlike for the case of the transverse train track coordinates, we
do not need to pass to a quotient space.
It is harder to see the inequalities satisfied by the shear coordinates, though.

When one train track maps to another, there is an induced linear
transformation of shear coordinates, so as with the transverse coordinates
we can take the union of the domains for which they work on any
particular train track.
Intersection number of a foliation defined by shear coordinates
with a lamination defined using tangential train track coordinates
is given by a simple formula, which turns out to be
the same formula that represents the symplectic pairing for that train track.
Abstractly, the transverse coordinates are in the dual space
of the tangential coordinates; the shear coordinate can be
defined by mapping the dual space to the tangential space using
the symplectic pairing.

One advantage of the shear coordinates is the uniqueness: they are determined
directly, and not just defined up to switch vectors.
The shear coordinates are particularly natural when $\mu$ is the one-skeleton
of an ideal triangulation.  They are also natural for describing earthquakes.
For any transverse invariant measure $\nu$ supported on $\mu$,
the earthquakes along $\nu$ are defined
by adding multiples of the transverse measure of $\nu$ to the initial
shear coordinates.  Just as for transverse coordinates,
the stretch maps are defined by multiplying the shear coordinates by scalars.
Another advantage of shear coordinates is that the formula for derivatives of
lengths of laminations along an earthquake mentioned (see
\ref{Convex models for measured lamination space}) generalizes:
for a close train track approximation $\tau$ of $\mu$, the directional
derivative of the length of a measured lamination $\lambda$ in the 
direction of a vector $V$ (in shear coordinates) is the integral
over the intersections of $\tau$ with $\lambda$ of the cosine
of the angle of intersection times the total shear along the branch of
$\tau$, integrated with respect to the transverse measure of $\lambda$.
This approximation can also be expressed directly in terms of 
$\mu$ and $\lambda$ --- the shear can be thought of as a transverse
invariant distribution for $\mu$, or a kind of current.
The formula is not hard to prove, by constructing an $\epsilon$-curved
embedding of $\lambda$ into the surface obtained by a small cataclysm along 
$\mu$.

We will not go into detail, but it is worth mentioning that
the complexified version of shear coordinates have an interpretation
as bending coordinates which parametrize a quasi-fuchsian 3-manifold
in terms of the 3-dimensional shape of its $\mu$ -pleated surface.

\section {The maximum stretch lamination is almost always a curve} \label {The maximum stretch lamination is almost always a curve}

As we have noticed, there are flat places on the boundary of the unit ball for
the stretch norm.  When a flat portion of the
boundary of a convex set has maximum possible dimension so that it has
interior, we will refer to it as a {\it facet}.

A facet determines a point in the dual unit sphere,
representing the linear function whose value is
1 on the facet.  Such points are recognized by the property that the dual
ball has a sharp point there --- that is, the tangent cone contains no lines.
The {\it Hausdorff dimension} $\delta ( X ) $ of a metric space $X$
is the infimum of the set of real numbers $d$ such that for all $\epsilon$
there is a countable cover of $X$ by balls of radius $r_{i} < \epsilon$
such that
$$ \sum  r_i^d < \epsilon . $$
The motivation for the definition of
Hausdorff dimension is that for a smooth manifold $M^{n}$ of dimension $n$, 
$r_{i}^{n}$ is roughly proportional to the volume of the $i $th ball, from
which one can easily deduce that the Hausdorff dimension is $n$.

Any subset of a manifold $M$ of Hausdorff dimension $d < n$
(with respect to the induced metric) has n-dimensional Lebesgue measure 0.
The difference $n-d$ is its {\it Hausdorff codimension},
which gives a good qualitative measure of how ``rare'' the set is.

\begin{theorem}[mostly facets]\label{mostly facets}
Almost every point on the
unit sphere of the stretch norm for the tangent space 
of Teichm\"uller space is on a facet.
In fact, the complement of the union of facets has Hausdorff codimension 1.

The facets correspond to simple closed geodesics.
In other words, each facet is contained in a plane 
$d \log \length ( \alpha ) = 1$ for some simple closed geodesic $\alpha$.
\end{theorem}

We will make some preparations before launching into the proof of the theorem.
A global version of this statement, which says that the maximally stretched
lamination for a pair of surfaces is almost always a simple closed curve,
will be given near the end of this section.

The first step of the proof will be to show that in terms
of their physical appearance on the surface,
laminations come only in very tightly grouped clusters.
From this, we will deduce that a random infinitesimal
change in a hyperbolic structure is unlikely to maximally
stretch anything except a simple closed curve.

Recall that the {\it Hausdorff metric}
measures the discrepancy between two closed subsets $A$ and $B$ of a metric
space as the supremum of the distance of a point in either set from the
closest point in the other set.

Let $\CL ( S )$ denote the set of maximal chain recurrent laminations
on $S$.  Their complementary regions are either ideal triangles or
once-punctured monogons.
The punctured monogons remain because a chain recurrent lamination can
never have a leaf going out to a cusp.

\begin{theorem}[CL has dimension 0]\label{CL has dimension 0}
For any hyperbolic surface of finite area $S$, the Hausdorff dimension of 
$\CL ( S )$ with the Hausdorff metric is 0.
\end{theorem}

Slightly more work would prove that the space of all geodesic laminations on 
$S$ have Hausdorff dimension 0.  This is closely related to a result of Joan
Birman and Caroline Series, which says that the
union of all geodesic laminations has measure 0, in fact Hausdorff dimension 1.
The two facts are related via a larger space, consisting of pairs 
$ ( \lambda , x )$ where $\lambda$ is a geodesic lamination and $x$ is a point
on $\lambda$.  The larger space can be shown to have Hausdorff dimension 1.
Since projection to the surface does not increase distance,
its image in the surface also has Hausdorff dimension 1.

\proof %{Proof of \shortlabel CL has dimension 0:10.2:}
Let $\tau$ be any maximal $\epsilon$-curved train track on $S$, and let 
$M$ be the space of invariant weights for $\tau$.

For any $R > 0$ , let $A_{R}$ be the set of homotopy classes of paths carried by$\tau$ which end in both directions at divergent switches.  For each 
$a \element A_{R}$, there is associated a linear function $\delta ( a , * )$
as in the proof of theorem \ref{Linear fragments from chain recurrence}
which measures the linear deviation of the holonomy from the case when there
is a leaf connecting the branch points for the associated highway.

The zero-sets of these linear functions divide the convex set $M$ into chambers.
A set of $m$ hyperplanes in convex set of a fixed dimension divides
it into a set of chambers whose cardinality is bounded by a polynomial
in the number of hyperplanes.  The number of elements of $A_{R}$, which define
the hyperplanes, is bounded by a polynomial in $R$.  By the proof of Theorem
\ref {Linear fragments from chain recurrence},
the Hausdorff distance between any two measured laminations which are also
in $\CL$ is not more than the order of $e^{-R}$.

For any real number $d > 0$, $N_{R} ( e^{-R} )^{d}$ goes to 0 as $R$
goes to infinity.  Since the elements of $\CL ( S )$ which are in the image of
$\ML ( S )$ are dense in $\CL$, this shows that the Hausdorff dimension of 
$\CL$ is 0.
\endproof %{\ref {CL has dimension 0}}

We will also need the generalization of this to more general
chain recurrent laminations, which are not necessarily maximal.
If $\tau$ is any recurrent, $\epsilon$-curved train track,
define $\CL_{\tau}$ to be the set of chain recurrent laminations
$\lambda$ carried on $\tau$ for which a regular neighborhood of
$\tau$ is a topologically stable neighborhood of $\lambda$.
The proof above gives

\begin{theorem}[CL tau has dimension 0]\label{CL tau has dimension 0}
For any $\tau$ as above, $\CL_{\tau}$ has Hausdorff dimension 0.
\end{theorem}

The way to look at this result which has significance for the proof of theorem
\ref {mostly facets} is that the Hausdorff dimension of 
$\CL_{\tau}$ is smaller than the dimension of the space of projective laminations
carried by $\tau$ in every case except the most trivial, when $\tau$ is a
simple closed curve.

For each $\lambda  \element  \CL ( S )$, there is associated a cataclysm
coordinate system.  When $S$ has cusps, the coordinate system is
defined by augmenting $\lambda$ with one new leaf going to each cusp.
Thus, there is a continuous map $c :  U\arrow \teich ( S )$ where 
$U  \subset  \CL ( S ) \cross \MF ( S)$ consists of pairs $( \lambda , F )$
such that $F$ is a measured foliation transverse to $\lambda$ and
standard near cusps, and where $\teich ( S )$ is the Teichm\"uller space of 
$S$.

To go to the infinitesimal level, we can first localize to a train track $\tau$.
Let $U_{\tau}$ consist of pairs $( \lambda , F )$ where $\lambda$ is an element
of $\CL ( S )$ carried by $\tau$ and $F$ is transverse to $\tau$.  The
$\MF$-factors within $U_{\tau}$ now have consistent linear structures,
so there is a vector bundle $E_{\tau}$ over $U_{\tau}$ consisting of the
tangent spaces to the $\MF$-factors.

Of course, $E_{\tau}$ does not really depend on $\tau$.  The $E_{\tau}$
piece together to give a globally defined vector bundle $E$ over $U$.

The infinitesimal version of the cataclysm parametrization $c$ is a map
$$ dc :  E\arrow T_* ( \teich ( S ) )
$$

\begin{theorem}[cataclysm parameters are infinitesimally Holder]\label{cataclysm parameters are infinitesimally Holder}
For any hyperbolic structure $g$ on the surface $S$, and any constant
$0 < \alpha < 1$, the map $dc$ defined above is $\alpha$-H\"older on the portion
of $E$ corresponding to hyperbolic structure $g$.
\end{theorem}

The Hausdorff metric on $\CL ( S )$ is defined in terms of a metric on $S$,
which is the reason for stating the theorem in terms of a particular hyperbolic
structure $g$.

\proof %{Proof of \shortlabel cataclysm parameters are infinitesimally Holder:10.4:}
The infinitesimal cataclysm map affects the developing map along a leaf $\gamma$
of the transverse measured foliation according to the derivative of the
linearly ordered infinite product formula of the proof of proposition
\ref {Neighborhood of a lamination}
The terms of the product correspond to intersections of 
$\gamma$ with ideal triangles in the complement of the maximal lamination
$\lambda$.

The tangent space of a Lie group at any point
is identified with the tangent space
at the origin by the action of the Lie algebra via
infinitesimal pre-composition (multiplication on the right)
in the Lie group.  Using this identification,
the Leibnitz rule for the derivative of a linearly ordered product
$$ \product_i g_i
$$
taken in an order $<$ in a nonabelian group becomes
$$ d \product_i g_i = \sum_i \left \{
\product_{ j > i } g_j \right \}^{-1} d g_i \product_{ j > i } g_j .
$$
(The notation here for a product involving both elements of the Lie
group and the Lie algebra is the limit of the macroscopic notation;
this is 
$$ \sum_j Ad \left ( \product_{ j > i } g_j \right ) d g_i
$$
in more standard notation).

How does the sequence $d g_{i}$ behave
as the measured foliation is changed infinitesimally? The norm $d g_{i}$
is proportional to the norm of $g_{i}$ times the derivative
of the transverse measure from $\gamma$ to the singular point
of the foliation within an ideal triangle.
This transverse measure is a linear function of the measured
foliation, so its derivative is bounded by a constant
times the norm of the element of the vector bundle $E$
times the transverse measure.  The linearly growing terms in the estimate for 
$d g_{i}$ are drowned by the exponential decresing terms $| g_{i} |$,
so the sequence $d g_{i}$ decreases exponentially.
The series for the derivative of the holonomy therefore converges
geometrically.

Consider two laminations $\lambda_{1}$ and $\lambda_{2}$ which are within 
$\epsilon$ of each other in the Hausdorff metric, and two tangent vectors 
$V_{1}$ and $V_{2}$ to $\MF$ which are also within $\epsilon$.  Any closed curve 
$\gamma$ on $S$ is homotopic to a stairstep curve $\gamma_{i} \ [i=1,2]$
which is made up of leaves of $\lambda_{i}$ and leaves of $F ( \lambda_{i} )$.
The two curves $\gamma_{i}$ can be chosen so that they stay
within a distance of each other which is bounded by a constant times
$\epsilon$ times the length of $\alpha$.  The holonomy about $\gamma$
can be computed in terms of the linearly ordered product using 
$\gamma_{i}$.  The products are taken over two different linearly ordered sets.
However, for each term of either product which is bigger than a constant times 
$\epsilon$, there is a corresponding term of the other product
which is within a constant times $\epsilon$.

Similarly, for the larger terms in the sum which expresses the infinitesimal
change in holonomy about $\alpha_{i}$ as $F ( \lambda_{i} )$ moves in the
direction $V_{i}$, there are corresponding nearby terms in the other sum.

We can compare the two series by looking separately at the large terms
and the small terms.  The terms of size less than a constant times 
$\epsilon$ sum to less than some other constant times $\epsilon$, since the
convergence is geometric.  There are the order of $\log ( \epsilon ) $
terms of size bigger than a constant times $\epsilon$, and each changes by at
most the order of $\epsilon$.  Therefore, the difference of the two
series has the order of $\epsilon \log ( \epsilon )$.

It follows that the derivative map $dc$ is H\"older of every constant less than 
$1$.
\endproof %{\ref {cataclysm parameters are infinitesimally Holder}}

Now we can show that at least the stretch vector fields themselves
take up a very small part of the tangent space to Teichm\"uller space.

\begin{theorem}[Stretch vector fields have dimension 0]\label{Stretch vector fields have dimension 0}
The set of directions in the tangent sphere bundle of Teichm\"uller space
which are tangent to stretch paths have Hausdorff dimension 0.
\end{theorem}

\proof %{Proof of \shortlabel Stretch vector fields have dimension 0:10.5:}
Consider the composition
$$ \CL ( S ) \arrow  E \arrow   T_1 ( \teich ( S ) ).
$$
The first arrow is associates to a lamination $\lambda$ in $\CL$ the pair 
$( \lambda ,  F ( \lambda ) )$ where $F ( \lambda )$ is
identified with a tangent vector to $\MF $.  The composition associates to 
$\lambda$ the tangent vector to its unit speed stretch path.

Clearly $F ( \lambda ) $ depends in a Lipschitz way on the Hausdorff
metric on $\CL$, since the foliation changes direction only as fast
as $\lambda$ changes direction.

The image of a set under a Lipschitz map has Hausdorff dimension
not greater than the Hausdorff dimension of its domain.
The image of a set under a $\alpha$-H\"older has dimension not greater than 
$1 / \alpha$ times the Hausdorff dimension of its domain, so $dc$
does not increase Hausdorff dimension.
Consequently, the image of the composition has Hausdorff dimension
not greater than 0.
\endproof %{\ref {Stretch vector fields have dimension 0}}

We will generalize this proposition to a statement involving laminations
which are not maximal.  For this, we localize to an $\epsilon$-curved train
track $\tau$.

We will define a consistent way to extend laminations in $\CL_{\tau}$
to maximal chain recurrent laminations on $S$.  First we augment $\tau$
to a recurrent train track $\sigma$ whose
complementary regions are triangles and punctured monogons.
We can do this by only attaching new branches to 
$\tau$.  As usual, some care is needed in the case that $\tau$ has orientable
components to make sure that there are new branches serving both as
exits and as entrances.  We may assume that $\sigma$ also has small curvature.

Any lamination $\lambda \element \CL_{\tau}$ can now be canonically extended
to a maximal chain recurrent lamination by adding one new
leaf for each branch of $\sigma - tau$.  Since $\lambda$
is a maximal lamination carried by $\tau$, there is
only one possible asymptotic behaviour for each new leaf.

This defines a Lipschitz embedding of $\CL_{\tau}$ in 
$\CL_{\sigma} \subset \CL$.  Let $\MF_{\sigma} $
be the space of measured foliations standard near the cusps of $S$
which are transverse to $\sigma$.  For each $\lambda \element \CL_{\tau}$,
there is a subspace $F_{\lambda} \subset \MF_{\sigma}$
consisting of those foliations whose restriction to a neighborhood
of $\lambda$ is a multiple of $F ( \lambda )$.
The measure class of the restriction map is a linear function,
so each of these subspaces have codimension $k$ equal to the dimension
of the space $\PF_{\tau}$.

\begin{proposition}[F sub lambda codimension k]\label{F sub lambda codimension k}
The Hausdorff codimension of the union of the subspaces
$F_{\lambda}$ is the dimension of $\PF_{\tau}$.
\end{proposition}

The map which assigns to an element of $\CL_{\tau}$
the measured foliation in a neighborhood of 
$\tau$ is Lipschitz, so its image has Hausdorff dimension 0.
The cone on this image has Hausdorff dimension 1.
The entire picture is formed by taking the product with a vector space,
so it has Hausdorff codimension $k$.
\endproof %{\ref {F sub lambda dimension k}}

\proof %{Proof of \shortlabel mostly facets:10.1:}
Define $E_{\tau} \subset  E$ to consist of those pairs $ ( \lambda ,  F )$
such that $\lambda$ is in $\CL_{\tau}$ and $F$ is in $F_{\lambda}$.  Clearly, 
$E_{\tau}$ is a set of Hausdorff $n-k$, where $n$ is the dimension of
Teichm\"uller space, since it is contained
in the product of a set of Hausdorff dimension $n-k$
with a set of Hausdorff dimension 0.

If $V$ is a tangent vector to Teichm\"uller space and $\lambda$
is a measured lamination which has maximal logarithmic derivative of length,
then the measure class of $F ( \lambda )$ has first derivative in the direction 
$V$ equal to some multiple of itself.  In other words, $V$ is represented as the
image by $dc$ of some element of $E_{\tau}$.

Since $dc$ is H\"older of every constant less than 1, it follows that
for any recurrent $\epsilon$-curved train track $\tau$,
the set of tangent vectors which maximally stretch some lamination in 
$\CL_{\tau}$ have Hausdorff codimension equal to the dimension of $\MF_{\tau}$.

The set of all tangent vectors which maximally stretch
some lamination other than a simple closed curve is thus
contained in a countable union of subspaces of Hausdorff
codimension at least 1, so the entire set has
Hausdorff codimension 1.
\endproof %{\ref {mostly facets}}

Here is the promised global version of the theorem:
\begin{theorem}[usually simple]\label{usually simple}
The maximal maximally stretched lamination $\mu ( g , h )$
for a pair of hyperbolic structures $g$ and $h$ is usually a simple closed
curve.  In fact, the set of pairs $ (g , h )$ such that $\mu ( g , h )$
is not a simple closed curve has Hausdorff codimension 1.
The same is true if either component is held constant.
\end{theorem}

\proof %{Proof of \shortlabel usually simple:10.7:} In the proof of the
infinitesimal version of this theorem, we made use of the fact that
for any chain recurrent lamination $\lambda$, there is a subspace 
$T_{\lambda}$ of the tangent space to Teichm\"uller space
at any point consisting of tangent vectors such that
with respect to some maximal lamination $\mu$ containing $\lambda$,
the cataclysm coordinates change as a multiple of $F ( \lambda )$
near $\lambda$.  The subspace $T_{\lambda}$ contains all tangent vectors which
maximally stretch (or contract) $\lambda$.

These subspaces piece together to form a foliation of Teichm\"uller space.  In 
$\mu$-cataclysm coordinates, in fact, the foliation is linear.
In particular, it is a real analytic foliation.

We can prove \ref {usually simple} in the same way as
\ref {mostly facets}
provided we can show that the cataclysm coordinates depend on the parameters.
We cannot expect the dependence on the lamination to be
H\"older of any fixed constant with respect to the Hausdorff metric for geodesic
laminations on any fixed hyperbolic surface.
Instead, we should form a composite metric taking into
account the geometry of many surfaces.

Let $K$ be any compact subset of Teichm\"uller space.  Then the
{\it composite Hausdorff metric} $d_{K}$ for geodesic laminations on $S$
is defined to be the supremum of the Hausdorff distances for realizations
of the laminations on surfaces of $S$.

\begin{proposition}[composite dimension]\label{composite dimension}
The Hausdorff dimension of the space of chain recurrent
geodesic laminations is 0 with respect $d_{K}$, for any compact $K$.
\end{proposition}

\proof [Proof of \ref{composite dimension}]
The previous proof still works.
In fact, the map from the Hausdorff metric coming from one
to the composite Hausdorff metric is H\"older of some constant,
so the image still has Hausdorff dimension 0.
\endproof %{\ref {composite dimension 0}}
\begin{proposition}[cataclysms Lipschitz]\label{cataclysms Lipschitz}
Let $J$ be a compact subset of $  \CL ( S ) \cross  \MF ( S )$
such that the lamination is transverse to the measured foliation.
Then there is a compact subset $K$ of Teichm\"uller space such that
the cataclysm parametrization
$$ \cat : J \arrow  \teich ( S ) $$
is Lipschitz with respect to the product of the metric 
$d_{K}$ and a linear metric.
\end{proposition}

\proof [Proof of \ref{cataclysms Lipschitz}]
We can assume that $J$ has the form of a product $J_{1} \cross J_{2}$,
since every compact subset is contained in a union of products.  In that case,
the measured foliations which occur have a natural consistent linear structure.
We may assume that $J_{2}$ is convex with respect to this linear structure.
With these assumptions, we take for $K$ the image of $J$ by $cat$.

We already know that when $\lambda $ is held fixed,
the cataclysm parameters are real analytic as a function of
the measured foliation.  The first derivative of the parametrization with
respect to the lamination is uniformly bounded if $\lambda$
ranges over a compact set of laminations.
Therefore, all we need to show is that for a fixed measured
foliation, the hyperbolic surface is a Lipschitz function
of the measured lamination.

For this, we can use the inverse cataclysm parametrization.
Given a lamination $\lambda $ and a hyperbolic structure $g$, the inverse
parametrization $F_{g} ( \lambda )$ gives
measured foliation which depends analytically on the hyperbolic structure.
The derivative of the cataclysm parametrization
in the forward direction is
uniformly bounded from being singular, so in the
backward direction the derivative is uniformly bounded.

Given a measured foliation $F$ and two laminations $\lambda_{1}$ and 
$\lambda_{2}$, we first look at the surface $\cat ( \lambda_{1} ,  F )$.
The lamination $\lambda_{2}$ has Hausdorff distance less than 
$d_{K} ( \lambda_{1} ,  \lambda_{2} )$ on this surface, so 
$F' = F_{ \cat ( \lambda_1 , F ) } (\lambda_2 )$ is within a constant times 
$\epsilon$ of $F$.  We have
$\cat ( \lambda_{1} , F ) = \cat ( \lambda_{2} , F ')$, so 
$\cat (\lambda_{2} , F )$ is not farther than a constant times $\epsilon$
away, by the uniform continuity of $\cat$ with respect to $F$.
\endproof %{\ref {cataclysms Lipschitz}}

\proof [Continuation, proof of \ref {usually simple}]
The proof is completed by an argument parallel to the infinitesimal version.
We have a family of parameter systems for Teichm\"uller space
such that the Hausdorff dimension
of the union is the same as the Hausdorff dimension of any particular one.

First we show that the set of pairs of hyperbolic structures $( g,h) $
having $\mu ( g,  h )$ contained in $\CL_{\tau}$ has Hausdorff codimension
equal to the dimension of $\PF_{\tau} $.
There is a linear projection of the space of measured foliations
transverse to a train track containing $\tau$ to the measure
class in the vicinity of $\tau$.  For $\mu ( g , h )$ to be contained in 
$\CL_{\tau}$, it is necessary (though not sufficient) that there be a lamination 
$\lambda  \element \CL$ such that in $\lambda$-cataclysm coordinates,
this image of $g$ in $\MF_{\tau}$ be a scalar multiple of the image of $h$.
The set of all elements of $\CL \cross  \MF  \cross  \MF$ satisfying
this necessary condition has Hausdorff codimension $k$
(measured using an appropriate composite Hausdorff metric on $\CL$)
equal to the dimension of $\PF_{\tau}$, so the image by the Lipschitz map
$\cat$ in Teichm\"uller space
cross Teichm\"uller space has the same codimension.

The statement when one of the factors is held fixed is proven in a similar way.
\endproof %{\ref {usually simple}}

It is interesting to contrast the shape of the unit stretch ball
in the tangent space with that in the cotangent space.
Because of the uniqueness of the topological maximally-stretched
lamination, two measured laminations which share a flat
in the dual unit ball must represent two different measures
on the same underlying lamination.
Any such flat can have dimension at most half the dimension
of the whole tangent space.
It has been shown, by Howie Masur, Bill Veech, and Steve Kerckhoff
(in independent analyses) that the set of all such laminations
has measure 0, in the piecewise linear metric on $\MF$.
This metric is Lipschitz-related to a linear metric on the
cotangent space of Teichm\"uller space, so the union of all
flats for the dual unit ball has measure 0.

The dual statements help to clarify the pictures of these two unit balls:
the ball consisting of all vectors which stretch no lamination
at a proportional speed greater than 1 (relative to its length)
is mostly made up of facets;
if it is set down on a table in a random orientation, it is almost
certain to land on a single point at which the surface of the ball is smooth.
On the other hand, the dual unit ball consisting of all 1-forms 
$d t \log \length ( \lambda )$ for some $t \le 1$ has almost no flats;
if it is set down on a table, it is almost certain to land
on a point at which the ball has a conical shape.

The fact that the maximally stretched lamination for a pair
of hyperbolic surfaces is almost always a simple closed
curve suggests that there may be an efficient computational
procedure to find $\mu ( g ,  h )$.
It would be especially worthwhile to extend this analysis
of lengths of laminations to three dimensions, to a similar analysis for the
lengths of laminations in quasi-Fuschsian groups.
Such an analysis might yield an efficient procedure to
explore the quasi-Fuchsian deformation space
(considered as a subspace of the character space for the
surface group), by giving
a reasonable estimate of how far one is from its boundary.
\bibliographystyle{halpha}
\bibliography{stretch}
\end{document}